\documentclass[12pt]{article}
\usepackage[psamsfonts]{amsfonts}
\usepackage{amsmath,amssymb,amsthm}
\usepackage[dvips]{graphicx}
\usepackage{eepic}
\def\UseSection{
        \numberwithin{equation}{section}
	\theoremstyle{plain}
        \newtheorem{theorem}    {Theorem}[section]
        \DefineTheorems 
}

\def\DefineTheorems{
	
	\newtheorem{lemma}      [theorem] {Lemma}
	
	\newtheorem{prop}       [theorem] {Proposition}
	
	\newtheorem{cor}        [theorem] {Corollary}

	\theoremstyle{definition}
	\newtheorem{defn}       [theorem] {Definition}
	\newtheorem{conj}       [theorem] {Conjecture}

	\theoremstyle{definition}

}

\newcommand{\bt}   {\begin{theorem}}
\newcommand{\et}   {\end  {theorem}}
\newcommand{\bl}   {\begin{lemma}}
\newcommand{\el}   {\end  {lemma}}
\newcommand{\bp}   {\begin{prop}}
\newcommand{\ep}   {\end  {prop}}
\newcommand{\bc}   {\begin{cor}}
\newcommand{\ec}   {\end  {cor}}
\newcommand{\bd}   {\begin{defn}}
\newcommand{\ed}   {\end  {defn}}

\newcommand{\ba}   {\begin{array}}
\newcommand{\ea}   {\end  {array}}
\newcommand{\be}   {\begin{enumerate}}
\newcommand{\ee}   {\end  {enumerate}}
\newcommand{\bi}   {\begin{itemize}}
\newcommand{\ei}   {\end  {itemize}}

\def\eq#1\en{\begin{equation}#1\end{equation}}  
\def\eqsplit#1\ensplit{
	\begin{equation}\begin{split}#1\end{split}\end{equation}
	}
\def\eqalign#1\enalign{
	\begin{align}#1\end{align}
	}
\def\eqmul#1\enmul{
	\begin{multline}#1\end{multline}
	}
\newcommand{\eqarrstar} {\begin{eqnarray*}} 
\newcommand{\enarrstar} {\end{eqnarray*}} 
\newcommand{\eqarray}   {\begin{eqnarray}} 
\newcommand{\enarray}   {\end{eqnarray}}

\newcommand{\lbeq}[1]  {\label{e:#1}}
\newcommand{\refeq}[1] {\eqref{e:#1}}    

\makeatletter
\newcommand{\labelcounter}[2]{{%
	\stepcounter{#1}
	\protected@write\@auxout{}%
	{\string\newlabel{#2}{{\csname the#1\endcsname}{\thepage}}}%
	{\ref{#2}}
	}}
\makeatother
\newcommand{\sss}   { \scriptscriptstyle } 

\newcommand{\Bbold} {{\mathbb B}}

\newcommand{\Ebold} {{\mathbb E}}

\newcommand{\Nbold} {{\mathbb N}}

\newcommand{\Pbold} {{\mathbb P}}
\newcommand{\Qbold} {{\mathbb Q}}
\newcommand{\Rbold} {{\mathbb R}}

\newcommand{\Zbold} {{\mathbb Z}}

\newcommand{\Acal}   {\mathcal{A}}

\newcommand{\Fcal}   {\mathcal{F}}

\newcommand{\Ocal}   {\mathcal{O}} 
\newcommand{\Pcal}   {\mathcal{P}}

\newcommand{\Rd}    {{ {\Rbold}^d}}
\newcommand{\Zd}    {{ {\Zbold}^d }}

\newcommand{\spose}[1] {{\hbox to 0pt{#1\hss}} }
\newcommand{\ltapprox} {\mathrel{\spose{\lower 3pt\hbox{$\mathchar"218$}}
 \raise 2.0pt\hbox{$\mathchar"13C$}}}
\newcommand{\gtapprox} {\mathrel{\spose{\lower 3pt\hbox{$\mathchar"218$}}
 \raise 2.0pt\hbox{$\mathchar"13E$}}}

\UseSection  
\renewcommand{\to}      {\rightarrow}

\newcounter{countC}  
\newcounter{countR}  
\newcommand{\Eb}{\Ebold}
\newcommand{\Pb}{\Pbold}
\newcommand{\R}{\Rbold}
\newcommand{\Z}{\Zbold}
\newcommand{\N}{\Nbold}
\newcommand{\cnctd}{\longleftrightarrow}
\newcommand{\conn}{\longrightarrow}

\newcommand{\nn}{\nonumber}

\newcommand{\smallsup}[1] {{\scriptscriptstyle{({#1}})}}
\newcommand{\sinfty}{{\!\sss \infty}}
\newcommand{\ICSBM}{infinite canonical super-Brownian motion }
\newcommand{\icsbm}{ICSBM }
\newcommand{\SP}{{\rm SP}}

\newcommand{\bC}{{\bf C}}
\newcommand{\lamb}{\lambda}
\newcommand{\lambc}{\lamb_{\rm c}}
\newcommand{\vep}{\varepsilon}
\newcommand{\mP}{{\mathbb P}}

\title  {
        Infinite canonical super-Brownian motion and scaling limits
        }

\author{
Remco van der Hofstad \footnote{Department of Mathematics and
Computer Science, Eindhoven University of Technology, P.O.\ Box
513, 5600 MB Eindhoven, The Netherlands. E-mail {\tt
rhofstad@win.tue.nl}}}

\begin{document}

\maketitle

    \begin{abstract}
    We construct a measure valued Markov process which we call {\it
    \ICSBM\!\!}, and which corresponds to the canonical measure
    of super-Brownian motion conditioned on non-extinction. Infinite
    canonical super-Brownian motion is a natural candidate for the
    scaling limit of various random branching objects on $\Z^d$
    when these objects are (a) critical; (b) mean-field and (c)
    infinite. We prove that \icsbm is the scaling limit of the
    spread-out oriented percolation incipient infinite cluster
    above 4 dimensions and of incipient infinite branching
    random walk in
    any dimension. We conjecture that it also arises as the
    scaling limit in various other models above the upper-critical
    dimension, such as the incipient infinite lattice tree above 8
    dimensions, the incipient infinite cluster for unoriented
    percolation, uniform spanning trees above 4 dimensions,
    and invasion percolation above 6 dimensions. This paper also serves
    as a survey of recent results linking super-Brownian to scaling
    limits in statistical mechanics.
    \end{abstract}


\section{Introduction}
\label{sec-intro}
Over the past years, it has become clear that super-Brownian motion (SBM)
arises as the scaling limit in various critical systems. Convergence
towards SBM can be expected when (a) the system is
critical; (b) the particles in the system are moving, and
undergo (critical) branching; and (c) the interaction in the system
is weak. In practice, requirement (c) means
that the system should be above the {\it upper critical dimension},
where the scaling ceases to depend on the dimension. Therefore, a Gaussian
limit can be expected, and this Gaussian limit is super-Brownian motion.
Examples where such results have been (partially) proved
are lattice trees above 8 dimensions \cite{DS97, DS98, Holm04},
the voter model above 2 dimensions \cite{CDP00, CDP99}, the contact process
above 4 dimensions \cite{DP99, HSa04c}, oriented percolation above 4
dimensions \cite{HS02a} and percolation above 6 dimensions \cite{HS00a, HS00b}.
Often, though not always, the proof of the convergence to SBM uses the lace
expansion.

Super-Brownian motion is the principle example of a measure-valued Markov process
in a similar way as Brownian motion is the principle example of a diffusion.
SBM has attracted considerable attention in the probability literature,
and has been described in detail in several recent books and major
reviews \cite{Daws93,Dynk94,Ethe00,LeGa99c,Perk02}.
The canonical measure of super-Brownian motion is described in
\cite{LeGa99c,Perk02}, and describes the structure of a single
continuum tree embedded into $\R^d$ where particles undergo critical
branching at any time scale, and move according to Brownian motions.
See also \cite{Slad02} for an excellent nontechnical introduction
describing scaling limits and their relations to super-processes.

There are two versions of convergence of a single geometric object
towards super-Brownian motion. In the first, the size of the critical object,
for example the critical percolation cluster, is fixed to be equal to $N$,
and the scaling limit as $N\rightarrow \infty$ is taken. In this case, the scaling
limit is {\it integrated super-Brownian excursion} (ISE), which is SBM
conditioned to have total mass 1. In the second, we investigate what happens when
these objects do not die out for a long time. In this case, the scaling
limit is the {\it canonical measure of super-Brownian motion}. We will describe this
limit in detail in Section \ref{sec-ICSBM} below, as an introduction
to super-Brownian motion for the non-specialist.

The canonical measure of super-Brownian motion is the scaling limit of a
single critical branching random walk which starts at the origin and survives
for some positive rescaled time. This connection will be made precise in
Sections \ref{sec-IIBRW} and \ref{sec-ICSBM}. Since critical branching processes
die out almost surely, also the canonical measure for SBM
dies out almost surely. In the above-mentioned examples, it is expected that
the critical structures live a {\it finite} amount of time almost surely,
and this illustrates why SBM can serve as the scaling limit for these
objects above the upper critical dimension.

Super-Brownian motion, however, cannot describe the scaling limit of critical systems where
the structures almost surely live an {\it infinite} amount of time.
Examples of such systems are invasion percolation and uniform spanning trees.
The aim of this paper is to study a critical super-process which survives with probability
one, and can be obtained as an appropriate limit of super-Brownian motion. It is natural to
expect that this object, which we call {\it \ICSBM}\!\!, serves as
the universal scaling limit of systems that are (a) infinite, though critical; (b) undergo
branching and motion; and (c) have weak interaction. There are different
versions of such infinite structures, namely {\it incipient infinite
structures}, which are obtained by appropriate limiting procedures in models where the
structures are a.s.\ finite, and {\it infinite} structures, where such a limit is not
necessary. An example for the former is the incipient infinite percolation cluster, an
example for the latter is a single tree in the uniform spanning forest.

This paper contains four main parts. In the first part (Section
\ref{sec-IIBRW}), we define incipient infinite branching random walk,
which is branching random walk conditioned
on non-extinction. This can be seen as a warm-up problem for the
construction of \ICSBM (\icsbm\!\!) in the second part (Section
\ref{sec-ICSBM}). We also discuss properties of \icsbm there. In the
third part (Section \ref{sec-iicop}), we prove that the incipient infinite cluster
in oriented percolation above 4+1 dimensions converges to \icsbm\!\!.
The latter result is based on the results obtained in \cite{HHS02a},
which is in turn based upon the convergence of finite-dimensional distributions
proved in \cite{HS02a}. In fact, it is fair to say that this paper is
inspired by these two papers. Finally, in the fourth part (Sections
\ref{sec-conjII}--\ref{sec-conjI}), we conjecture that many other
models also scale to \icsbm\!\!, making \icsbm a universal object.

This paper has two aims. The first aim is to introduce and investigate
\ICSBM\!\!, to state results concerning convergence towards \icsbm\!\!,
and to make conjectures in models where we cannot (yet?) prove such
convergence. The second aim of this paper is to review super-Brownian
motion and the recent results on super-Brownian motion arising as
the scaling limit in various critical high-dimensional models. This paper will
attempt not to be technical, and is aimed for the non-specialists both in the field
of the applications, as well as in the field of super-processes.
As a warm-up, we start by defining incipient infinite branching random
walk, which is branching random walk conditioned on non-extinction, and
we will see that \icsbm serves as the scaling limit of incipient
infinite branching random walk in a similar way as super-Brownian
motion serves as a scaling limit for ordinary branching random walk.

\section{Incipient infinite branching random walk}
\label{sec-IIBRW}
In this section, we will construct the incipient infinite branching
random walk (IIBRW) measure, which is the measure of branching
random walk conditioned on non-extinction. Since SBM is the
scaling limit of branching random walk, it is instructive to perform
the construction first for branching random walk. This construction
is simpler due to the discrete nature of branching random walk, and its
close connection to branching processes. We can think of IIBRW as an
embedding of a critical branching process conditioned on
non-extinction into $\Z^d$. Branching process conditioned
on non-extinction have a long history, which we will review
in some detail below.

We have two constructions for IIBRW. We will see that the two
definitions of IIBRW coincide. In the next section, we will construct
\icsbm in two ways, and these constructions mirror the two constructions
for IIBRW given here.

\subsection{Model and main results}
We start by introducing branching random walk. We follow the construction
in \cite{BCHS99}. Branching random walk is defined in terms of embeddings
of abstract trees into $\Zd$. The abstract trees are the family trees of
the critical branching process with a critical offspring
distribution $(p_m)_{m=0}^{\infty}$ with finite variance.
For simplicity, we will assume that
$(p_m)_{m=0}^{\infty}$ has a finite third moment.

In more detail, we begin with a single individual having
$\xi$ offspring, where $\xi$ is a random variable with
distribution $(p_m)_{m=0}^{\infty}$, {\it i.e.},
$\mathbb{P}(\xi = m) = p_m$ with
    \eq
    \sum_m mp_m=1, \qquad \sigma_p^2=\sum_m m(m-1) p_m<\infty.
    \en
Each of the offspring
then independently has offspring of its own, with the
same critical distribution $(p_m)_{m=0}^{\infty}$.
For a tree $T$, with the $i^{\rm th}$ individual
having $\xi_i$ offspring, this associates to $T$ the
probability
    \eq
    \lbeq{PTdef}
    \mathbb{P}(T) = \prod_{i \in T} p_{\xi_i}.
    \en
The product is over the vertices of $T$.

It is important to be clear about when two trees $T$ are the same and
when they are not.  For this, we introduce a description of $T$ in
terms of {\em words}.  These words arise inductively as follows.
The root is the word $0$.  The children of the root are
the words $01, 02, \ldots ,0\xi_0$.
The children of $01$ are the words $011, \ldots , 01\xi_{01}$, and so on.
The family tree is then uniquely represented by a set of words.
Two trees are the same if and only if they are represented by the same
set of words.

We define an embedding $\phi$ of $T$ into $\Zd$ to be a mapping
from the vertices of $T$ into $\Zd$ such
that the root is mapped to the origin and, given that $i$ is mapped to $x\in \Zd$,
the child $j$ of $i$ is mapped to $y\in \Zd$ with probability $D(y-x)$.
We will always assume
that $D$ is symmetric, and that $D$ has finite variance, i.e.,
    \eq
    \lbeq{sigdef}
    \sigma^2 = \sum_{x \in \Zd} |x|^2 D(x)<\infty,
    \en
where $|\cdot |$ denotes the Euclidean norm on $\Rd$.
We will also assume that for some $\delta\in (0,1)$,
    \eq
    \lbeq{deltadef}
    \sigma^2 = \sum_{x \in \Zd} |x|^{2+2\delta} D(x)<\infty.
    \en
In later sections,
we will put stronger conditions on $D$, but for branching random walk,
this is not necessary.

Branching random walk is then defined to be the
set of configurations $(T, \phi)$, with probabilities
    \eq
    \lbeq{PTphidef}
    \mathbb{P}^{\smallsup{\rm brw}}(T, \phi) = \mathbb {P}(T) \prod_{ij\in T} D(\phi(j)-\phi(i)).
    \en
Here $ij\in T$ means that $j$ is the child of $i$ in the tree $T$.
In particular, the path in $\Zd$ from the origin to $\phi(i)$, where
$i\in T$ is a random walk path of length $|i|$ with transition probabilities
given by $D$. Here $|i|$ denotes the generation of $i$ in $T$, which
is the same as the graph distance between the root of $T$ and $i$.

Critical branching processes die out, i.e., $\mathbb {P}(|T|<\infty)=1$.
We now give two ways of generating a measure on infinite trees.
We let $T_m$ denote the restriction of $T$ to the points that
are at most tree distance $m$ away from the root, i.e.,
$T_m=\{i\in T: |i|\leq m\}$. Then, we let $(T,\phi)_m$
denote the embedding of $T_m$ in $\Z^d$.

Firstly, for a realisation $C$ of the embedded tree up to
time $m$, we define
    \eq\lbeq{PnBRWdef}
    \mathbb {P}_n^{\smallsup{\rm brw}}(C) =
    \sum_{x\in \Z^d}\sum_{i\in T:|i|=n} \mathbb {P}^{\smallsup{\rm brw}}((T,\phi)_m=C, \phi(i)=x),
    \en
and we let
    \eq
    \lbeq{BRWIICdef}
    \mathbb {P}_{\sinfty}^{\smallsup{\rm brw}}(C) =
    \lim_{n\rightarrow \infty} \mathbb {P}_n^{\smallsup{\rm brw}}(C),
    \en
assuming the limit exists. Denote by
    \eq
    \lbeq{Nndef}
    N_n=\#\{i\in T: |i|=n\}
    \en
the number of particles alive at time $n$. Then
    \eq
    \lbeq{altPnform}
    \mathbb {P}_n^{\smallsup{\rm brw}}(C) =
    \Ebold^{\smallsup{\rm brw}}(N_n I[(T,\phi)_m=C]),
    \en
where, for an event $E$, $I[E]$ denotes the indicator of $E$.
The measure $\mathbb {P}_n^{\smallsup{\rm brw}}$ is a probability measure, since
$\{N_n\}_{n=0}^{\infty}$ is a {\it martingale} with $N_0=1$. Therefore,
$\mathbb {P}_n^{\smallsup{\rm brw}}$ is a martingale transformation of
$\mathbb {P}^{\smallsup{\rm brw}}$.

For our second construction, we let
    \eq\lbeq{QnBRWdef}
    \mathbb {Q}_n^{\smallsup{\rm brw}}(C) =
    \mathbb {P}^{\smallsup{\rm brw}}((T,\phi)_m=C | \exists i\in T: |i|=n),
    \en
and we let
    \eq
    \lbeq{BRWIICaltdef}
    \mathbb {Q}_{\sinfty}^{\smallsup{\rm brw}}(C) =
    \lim_{n\rightarrow \infty} \mathbb {Q}_n^{\smallsup{\rm brw}}(C),
    \en
assuming the limit exists. The conditioning that there is an $i\in T$ with $|i|=n$
means that we condition the branching process to be alive at time $n$.

\begin{theorem}
\label{thm-BRWeq}
The measures in \refeq{BRWIICdef} and \refeq{BRWIICaltdef} are well-defined and
$\mathbb {P}_{\sinfty}^{\smallsup{\rm brw}}=\mathbb {Q}_{\sinfty}^{\smallsup{\rm brw}}$.
\end{theorem}

\proof 
The events $(T,\phi)_m=C$ are cylinder events, and we will prove convergence for these
cylinder events first. When $(T,\phi)_m=C$, we have that $N_m=N_m(C)$, which is the number of particles
in generation $m$ for $C$. We will first show
    \eq\lbeq{Pinftyeq}
    \mathbb {P}_{\sinfty}^{\smallsup{\rm brw}} (C)=N_m(C) \mathbb {P}^{\smallsup{\rm brw}}((T,\phi)_m=C)=
    \mathbb {P}_m^{\smallsup{\rm brw}}(C).
    \en
Equation \refeq{Pinftyeq} shows that $\mathbb {P}_{\sinfty}^{\smallsup{\rm brw}}$ can be
seen as a martingale change of measure of $\mathbb {P}^{\smallsup{\rm brw}}$.
Such transformations occur more generally for conditioned stochastic
processes, and are often called {\it $h$-transforms}. See \cite{SV00} for an example
where the $h$-transform is used to compute the super-Brownian motion exit measure.

We compute for every $n\geq m$
    \eqarray
    \lbeq{IIBRWid}
    \mathbb {P}_n^{\smallsup{\rm brw}}(C) &=& \mathbb {E}^{\smallsup{\rm brw}}(N_n I[(T,\phi)_m=C])\\
            &=& \mathbb {P}^{\smallsup{\rm brw}}((T,\phi)_m=C)\mathbb {E}^{\smallsup{\rm brw}}
                (N_n|(T,\phi)_m=C)\nonumber\\
            &=& \mathbb {P}^{\smallsup{\rm brw}}((T,\phi)_m=C)
            \mathbb {E}^{\smallsup{\rm brw}}(N_n|N_m=N_m(C))=N_m(C)
            \mathbb {P}^{\smallsup{\rm brw}}((T,\phi)_m=C),\nonumber
    \enarray
since $\{N_m\}_{m=0}^{\infty}$ is a martingale. As the right-hand side does not
depend on $n$, we also must have that \refeq{Pinftyeq} holds.

We will next show that also $\mathbb {Q}_{\sinfty}^{\smallsup{\rm brw}} (C)$
equals the right-hand side of \refeq{Pinftyeq},
which will prove Theorem \ref{thm-BRWeq}. We first let
    \eq
    \theta_n=\mathbb {P}^{\smallsup{\rm brw}}(\exists i\in T: |i|=n).
    \en
We again compute, for $n\geq m$, and using that $(N_n)_{n=0}^{\infty}$ is a martingale,
    \eqarray
    \mathbb {Q}_n^{\smallsup{\rm brw}}(C) &=& \frac{1}{\theta_n}\mathbb {P}^{\smallsup{\rm brw}}((T,\phi)_m=C, \exists i\in T:|i|=n)\lbeq{Qn}\\
            &=& \frac{1}{\theta_n}
            \mathbb {P}^{\smallsup{\rm brw}}((T,\phi)_m=C) \mathbb {P}^{\smallsup{\rm brw}}
            (\exists i\in T:|i|=n|(T,\phi)_m=C)\nonumber\\
            &=& \frac{1}{\theta_n} \mathbb {P}^{\smallsup{\rm brw}}((T,\phi)_m=C)
            \mathbb {P}^{\smallsup{\rm brw}}
            (\exists i\in T:|i|=n|N_m=N_m(C))\nonumber\\
            &=& \mathbb {P}^{\smallsup{\rm brw}}((T,\phi)_m=C)
            \frac{1-(1-\theta_{n-m})^{N_m(C)}}{\theta_n},\nonumber
    \enarray
where in the final equality, we use that the particles in the first generation
evolve independently. As $n\rightarrow \infty$, the right-hand side of \refeq{Qn} converges to
$N_m(C)\mathbb {P}^{\smallsup{\rm brw}}((T,\phi)_m=C)$
whenever $\theta_n\rightarrow 0$, and
$\frac{\theta_n}{\theta_{n-m}}\rightarrow 1$. In fact,
for branching random walk, we know a lot more (see e.g.
\cite{Aldo93a}), namely that
    \eq\lbeq{qnasy}
    \lim_{n\rightarrow \infty} n\theta_n = \frac{2}{\sigma_p^2}.
    \en
We now complete the proof of Theorem \ref{thm-BRWeq}. Cylinder events of the form
$(T,\phi)_m=C$ generate the $\sigma$-algebra of all events. Since the
limiting measure $\Pbold_{\sinfty}^{\smallsup{\rm brw}}$ is consistent, we can extend it
to the full $\sigma$-algebra by Kolmogorov's Extension Theorem (see e.g.\ \cite{Simo79}).
This completes the proof.
\qed


\subsection{The branching random walk higher-point functions}
A convenient way to describe the distribution of a discrete random measure is
by using the {\it $r$-point functions}. In this section, we will assume that
$(p_m)_{m=0}^{\infty}$ has all moments. The $r$-point function describe the
numbers and locations of particles present at various times.
Denote by
    \eq
    \lbeq{taurbrwdef}
    \tau_{n_1,\ldots,n_{r-1}}(x_1,\ldots, x_{r-1})
    = \sum_{i_1, \ldots, i_{r-1}}
    {\Pbold}^{\smallsup{\rm brw}}(i_j\in T, \phi(i_j)=x_j, |i_j|=n_j
    \mbox{ for each }
    j=1,\ldots,r-1)
    \en
the branching random walk $r$-point functions. We will often abbreviate
    \eq
    \tau_{\vec n}(\vec x)=
    \tau_{n_1,\ldots,n_{r-1}}(x_1,\ldots, x_{r-1}).
    \en
The $r$-point functions give rise to measures, and these measure are
called {\it mean moment measures}, or just {\it moment measures}.
Indeed, let $n\geq 0$, and define the random measures
$\{\mu_n\}_{n=0}^{\infty}$ by
    \eq
    \lbeq{mundef}
    \mu_n(x)=\sum_{i\in T: |i|=n} I[\phi(i)=x].
    \en
The random measures $\{\mu_n\}_{n=0}^{\infty}$ evolve in time and $\mu_n$
describes the amount of mass and the spatial location of the mass
of the BRW at time $n$. We start with a single particle at time
0 located at the origin, so that
    \eq
    \lbeq{mu0def}
    \mu_0(x)=\delta_{x,0}.
    \en
Then, the law of the measured-valued process $\{\mu_n\}_{n=0}^{\infty}$
can be computed in terms of
the joint moments $\Ebold^{\smallsup{\rm brw}}[\prod_{i=1}^{r-1} \mu_{m_i}(y_i)^{a_i}],$
where $m_i\in \N, y_i\in \Z^d, a_i\in \N$. The $r$-point functions appear
explicitly in this description, since
    \eq
    \Ebold^{\smallsup{\rm brw}}[\prod_{i=1}^{r-1} \mu_{m_i}(y_i)^{a_i}]
    =\tau_{\vec n}(\vec x),
    \en
where $(x_j,n_j)$ equals $(y_i, m_i)$ precisely $a_i$ times. Thus, the joint moments
of the measures $\{\mu_n\}_{n=0}^{\infty}$ are equal to the $r$-point functions.
In the remainder of this section, we will give recursive formulas for $\tau_{\vec n}(\vec x)$,
and identify the $r$-point functions of IIBRW in terms of those of BRW.
We start with the latter.
The IIBRW $r$-point functions are defined, for $m_i \geq 0$ and $x_i \in \Zd$, and similarly
to \refeq{taurbrwdef}, by
    \eq
    \lbeq{mpoint,nobranchesbrw}
    \rho_{\vec{m}}(\vec{x})
    = \sum_{i_1, \ldots, i_{r-1}}
    {\Pbold}_{\sinfty}^{\smallsup{\rm brw}}(i_j\in T, \phi(i_j)=x_j, |i_j|=m_j
    \mbox{ for each }
    j=1,\ldots,r-1).
    \en
In the theorem below we identify $\rho_{\vec{m}}(\vec{x})=\rho_{m_1,\ldots,m_{r-1}}(x_1,\ldots, x_{r-1})$:

\begin{theorem}
\label{prop-relrhotaubrw}
For all $\vec{m}=(m_1, \ldots, m_{r-1})$ with $m_i\geq 0$ and $\vec{x}=
(x_1, \ldots, x_{r-1})$ with $x_i\in \R^d$,
    \eq
    \lbeq{relrhotaubrw}
    \rho_{m_1,\ldots,m_{r-1}}(x_1,\ldots, x_{r-1})
    =\sum_{x_0\in \Z^d}\tau_{\bar{m}, m_1,\ldots,m_{r-1}}
    (x_0,x_1,\ldots, x_{r-1})
    ,
    \en
where $\bar{m}$ denotes the largest component of $\vec{m}=(m_1,\ldots,m_{r-1})$.
\end{theorem}


\proof
To prove \refeq{relrhotaubrw} for branching random walk, we observe
that by \refeq{taurbrwdef} and \refeq{Pinftyeq} and the fact that
the event that $\{i_j\in T, \phi(i_j)=x_j,
|i_j|=m_j \text{ for each }j=1, \ldots, m-1\}$ only depends on $(T,\phi)_m$,
        \eqalign
        &\rho_{m_1,\ldots,m_{r-1}}(x_1,\ldots, x_{r-1})\\
        &=\lim_{n \rightarrow \infty}
        \sum_{i_1, \ldots, i_{r-1}}
        {\Pbold}_n^{\smallsup{\rm brw}}(i_j\in T, \phi(i_j)=x_j, |i_j|=m_j
        \mbox{ for each }
        j=1,\ldots,r-1)
        \nn\\
        &=\lim_{n \rightarrow \infty}
        \sum_{i_1, \ldots, i_{r-1}}
        {\Pbold}_{\bar m}^{\smallsup{\rm brw}}(i_j\in T, \phi(i_j)=x_j, |i_j|=m_j
        \mbox{ for each }
        j=0,\ldots,r-1)
        \nn\\
        &=
        \sum_{x_0} \tau_{\bar{m},\vec{m}}(x_0, x_1, \ldots, x_{r-1}),\nn
        \enalign
where in the second equality, we write $m_0=\bar{m}$.

%
    \qed


By Theorem \ref{prop-relrhotaubrw}, to identify the $r$-point functions
of IIBRW, it suffices to identify the $r$-point functions of BRW. We
will now investigate the $r$-point functions for BRW.
We first introduce some notation. Let $(f_j)_{j=0}^{\infty}$ denote the factorial
moments of the distribution $(p_m)_{m=0}^{\infty}$, i.e.,
    \eq
    f_j = \sum_{m=j}^{\infty} \frac{m!}{(m-j)!} p_m.
    \en
Also, we write $\Pcal_j$ for the number of partitions of
$\{1, \ldots, r-1\}$ into $j$ non-empty sets, where we order the elements
of $\vec I\in \Pcal_j$ by ordering the smallest components. Thus, $I_1$ contains
the element 1. Finally, for $I=\{i_1, \ldots, i_j\}\subseteq \{1,\ldots, r-1\}$, we write
$\vec n_I=(n_{i_1}, \ldots, n_{i_j})$. We will prove the following proposition:

\begin{prop}
\label{prop-taur}
For every $\vec{x}\in \Z^{d(r-1)}$ and every $\vec{n}=(n_1, \ldots, n_{r-1})$ with
$n_i\geq 1$ for all $i=1, \ldots, r-1$,
    \eq
    \lbeq{taurec}
    \tau_{\vec n}(\vec x)=\sum_{j=1}^{r-1} f_j \sum_{\vec{I}\in \Pcal_j}
    \prod_{s=1}^{j} (D*\tau_{\vec n_{I_s}-1})(\vec{x}_{I_s}).
    \en
\end{prop}
Before proving Proposition \ref{prop-taur}, we will discuss its relevance.
The significance of \refeq{taurec} lies in the fact that we can use it recursively
to identify the $r$-point functions. As a side remark, we immediately see that when
$f_{r-1}=\infty$, then also there exist $x_1, \ldots, x_{r-1}$ such that
$\tau_{\vec n}(\vec x)=\infty$. For
$r=2$, we obtain
    \eq
    \lbeq{tauneq}
    \tau_{n}(x)=f_1 (D*\tau_{n-1})(x), \qquad \text{ so that }
    \qquad \tau_n(x) = f_1^n D^{*n}(x)=D^{*n}(x),
    \en
where the last equality holds since the branching process is critical.

A special example arises when we consider binary branching, i.e., $p_m=
\frac 12 (\delta_{m,0}+\delta_{m,2})$. In this case, $f_1=f_2=1,$ and
$f_m=0$ for all $m\geq 3$.\footnote{The fact that
$f_m=0$ for all $m\geq 3$ is implied by $p_m=0$ for $m\geq 3$, and thus, the
computation holds somewhat more generally.}
Thus, we obtain that, writing $I=I_2$, so that $1\not\in I$,
    \eq
    \lbeq{taurec2}
    \tau_{\vec n}(\vec x)=(D*\tau_{\vec n-1})(\vec x)+\sum_{I\subseteq J_1: I\neq \varnothing}
    (D*\tau_{\vec n_{I}-1})(\vec{x}_{I})(D*\tau_{\vec n_{J\backslash I}-1})(\vec{x}_{J\backslash I}),
    \en
where $J=\{1, \ldots, r-1\}, J_1=J\backslash \{1\}.$
Iterating the recursion yields
    \eq
    \lbeq{taureq2}
    \tau_{\vec n}(\vec x)=\sum_{I\subseteq J_1: I\neq \varnothing} \sum_{m=0}^{\underline{n}-1} \sum_y D^{*m}(y)
    (D*\tau_{\vec n_{I}-m-1})(\vec{x}_{I}-y)(D*\tau_{\vec n_{J\backslash I}-m-1})(\vec{x}_{J\backslash I}-y),
    \en
where $\underline{n}$ denotes the minimal element of
$\vec n=(n_1, \ldots, n_{r-1})$.
Using \refeq{tauneq}, we can write \refeq{taureq} as
    \eq
    \lbeq{taureq}
    \tau_{\vec n}(\vec x)=\sum_{I\subseteq J_1: I\neq \varnothing} \sum_{m=0}^{\underline{n}-1} \sum_y \tau_m(y)
    (D*\tau_{\vec n_{I}-m-1})(\vec{x}_{I}-y)(D*\tau_{\vec n_{J\backslash I}-m-1})(\vec{x}_{J\backslash I}-y).
    \en
Equation \refeq{taureq} yields an explicit recursion for the $r$-point function in
terms of $r$, since on the right-hand side only $s$-point functions with
$s<r$ appear. For different offspring distributions, \refeq{taurec} is not
so easily solved, and in Section \ref{sec-MMSBM} below, we will identify
the scaling limit of $\tau_{\vec n}(\vec x)$ for general offspring distributions, by proving that
the contribution due to $j\geq 3$ in \refeq{taurec} is an error term.
\vskip 0.5cm

\noindent
{\it Proof of Theorem \ref{prop-taur}.} Recall that
    \eq
    \lbeq{taurbrwrec}
    \tau_{n_1,\ldots,n_{r-1}}(x_1,\ldots, x_{r-1})
    = \sum_{i_1, \ldots, i_{r-1}}
    {\Pbold}^{\smallsup{\rm brw}}(i_j\in T, \phi(i_j)=x_j, |i_j|=n_j
    \mbox{ for each }
    j=1,\ldots,r-1)
    \en
We fix $i_1, \ldots, i_{r-1}$ in \refeq{taurbrwrec}. We condition on the number of offspring
of the root, and denote this number by $l$. These particles
are labeled as $01, 02, \ldots, 0l$. We write, for $i=1, \ldots, l$,
    \eq
    A_i=\{j: 0i\conn i_j\},
    \en
to be the indices that $0i$ is connected to. Thus, $j\in A_i$ precisely when $0i$ is an ancestor
of $i_j$. We have that $A_i\cap A_j=\varnothing$ for $i\neq j$
and $\bigcup_{i=1}^l A_i=\{1, \ldots, r-1\}$. Denote by $y_i$ the spatial location of
$0i$. Then we can write, for each $\Acal_1, \ldots, \Acal_l$,
    \eqalign
    &\sum_{i_1, \ldots, i_{r-1}}
    {\Pbold}^{\smallsup{\rm brw}}(i_j\in T, \phi(i_j)=x_j, |i_j|=n_j
    \forall
    j=1,\ldots,r-1,  A_t=\Acal_t \forall t=1,\ldots, l|\xi_0=l, \phi(0t)=y_t)\nn\\
    &\qquad=\prod_{t=1}^l \sum_{i_j: j\in \Acal_t}
    {\Pbold}^{\smallsup{\rm brw}}(i_j\in T, \phi(i_j)=x_j-y_t, |i_j|=n_j-1
    \mbox{ for each }
    j\in \Acal_t)\nn\\
    &\qquad= \prod_{t=1}^l \tau_{\vec{n}_{\Acal_t}-1}(\vec{x}_{\Acal_t}-y_t).
    \enalign
Then, we end up with
    \eqalign
    \tau_{\vec n}(\vec x)
    &=\sum_{l=1}^{\infty}\sum_{\vec{A}}
    \sum_{y_1, \ldots, y_l}
    \Big[\prod_{t=1}^l \tau_{\vec{n}_{A_t}-1}(\vec{x}_{A_t}-y_t)\Big]
    {\Pbold}^{\smallsup{\rm brw}}(\xi_0=l, \phi(0t)=y_t)\nn\\
    &=\sum_{l=1}^{\infty}\sum_{\vec{A}}
    p_l \sum_{y_1, \ldots, y_l}\prod_{t=1}^l D(y_t)\tau_{\vec{n}_{A_t}-1}(\vec{x}_{A_t}-y_t)
    =\sum_{l=1}^{\infty}\sum_{\vec{A}} p_l
    \prod_{t=1}^l (D*\tau_{\vec{n}_{A_t}-1})(\vec{x}_{A_t}).
    \enalign
We further note that when $A_i=\varnothing$, then $\tau_{\vec{n}_{A_t}}(\vec{n}_{A_t}-y_t)
=1$, so that we can restrict the product over $A_t$ such that $A_t\neq\varnothing$.
Let $j$ denote the number of non-empty elements of $A_t$, and write
$I_1, \ldots, I_j$ for the non-empty elements of $A_t$, ordered in the unique way that
$\vec I\in \Pcal_j$. We can identify $I_s=A_{t_s},$
so that we have
    \eq
    \prod_{t=1}^l (D*\tau_{\vec{n}_{A_t}-1})(\vec{x}_{A_t})
    =\prod_{s=1}^j (D*\tau_{\vec{n}_{A_{t_s}}-1})(\vec{x}_{A_{t_s}})
    =\prod_{s=1}^j (D*\tau_{\vec{n}_{I_s}-1})(\vec{x}_{I_s}).
    \en
Then, the number of different ways of choosing
$A_1, \ldots, A_l$ such that $I_1, \ldots, I_j$ are fixed equals
$\frac{l!}{(l-j)!}$ for each $\vec{I}=(I_1, \ldots, I_j)\in \Pcal_j$. Thus, we arrive at
    \eqalign
    \tau_{\vec n}(\vec x)
    &=\sum_{j=1}^{r-1} \sum_{\vec{I}\in \Pcal_j} \sum_{l=1}^{\infty} p_l \frac{l!}{(l-j)!}
    \prod_{s=1}^j (D*\tau_{\vec{n}_{I_s}-1})(\vec{x}_{I_s})
    =\sum_{j=1}^{r-1} \sum_{\vec{I}\in \Pcal_j} f_j
    \prod_{s=1}^j (D*\tau_{\vec{n}_{I_s}-1})(\vec{x}_{I_s}).
    \enalign
This completes the proof of Proposition \ref{prop-taur}.
\qed

\subsection{The immortal particle}
\label{sec-BRWIP}
In the following theorem, we investigate the number of particles at any given time
that have infinitely many offspring. In its statement, and for $i\in T$, we write $i\conn \infty$
for the statement that the tree $T(i)$ rooted at $i$ is infinite.
    \begin{theorem}
    \label{thm-immparti}
    Under $\mathbb {P}_{\sinfty}^{\smallsup{\rm brw}}$, for every $m$, there is a unique
    $i\in T$ with $|i|=m$ such that $i\conn \infty$.
    \end{theorem}

The above result says that there is a single {\it immortal
particle}. This immortal particle performs a random walk
with transition probabilities $D$, and the mass alive at
any time is produced along the path of this single particle
and performs unconditioned critical branching random walk.
Before proving Theorem \ref{thm-immparti}, we first explain
this immortal particle picture in detail.

We construct IIBRW in the following way. Let $V_0$
be the root of the tree, and let $V_n$ be the label of the (unique)
child of $V_0V_1\cdots V_{n-1}$ that has an infinite tree emerging
from it. Let $\zeta_n$ be the total offspring of $V_{n-1}$. Then,
$\{(V_n,\zeta_n)\}_{n=0}^\infty$ is an i.i.d.\ sequence
with law
    \eq
    \lbeq{Vzetalaw}
    \Pbold^{\smallsup{\rm brw}}_{\sinfty}(V_n=j,\zeta_n=k)= p_{k+1} \qquad (1\leq j\leq k+1).
    \en
The (unique) infinite line of decent is now
$V_0, V_0V_1, V_0V_1V_2, \ldots.$ Embed this infinite paths into
$\Z^d$ as a single random walk path. Then, at the node $V_0V_1\cdots V_n$
in the infinite path, start $\zeta_n$ independent unconditioned branching
random walks, for which the root has word $V_0V_1\cdots V_{n-1} j$ for any
$j\neq V_n$. The law of the obtained process is equal to the law of
$(T,\phi)$ under $\mathbb {P}_{\sinfty}^{\smallsup{\rm brw}}$.

The above construction is quite involved, as we need to keep track of
what the infinite line of decent is. If we were only to be interested
in the spatial locations of the particles $(\phi(i))_{i\in T}$
rather than in the tree together with the spatial locations of the
particles $(T,\phi)$, then the construction simplifies considerably.
Indeed, in this case, we create single infinite random walk path
$\{\omega(n)\}_{n=0}^{\infty}$ in $\Z^d$ and at each
position $\omega(n)$, we start $\zeta_n$ independent
unconditioned branching random walks, where $\{\zeta_n\}_{n=0}^\infty$
is an i.i.d.\ sequence with law
    \eq
    \lbeq{zetalaw}
    \Pbold^{\smallsup{\rm brw}}_{\sinfty}(\zeta_n=k)= (k+1)p_{k+1}\qquad (k\geq 0).
    \en
Therefore, $\zeta_n$ is a size-biased version of the law $(p_m)_{m=0}^{\infty}$ minus one,
and is the marginal of $\zeta_n$ in the law in \refeq{Vzetalaw}.
We will discuss the history of this problem in Section \ref{sec-discBRW} below.
\vskip0.2cm

\noindent
{\it Proof of Theorem \ref{thm-immparti}.}
Fix $k\geq m$. We write $i\conn n$ for the event that there exists
$j\in T(i)$ with $|j|=n$. Then,
    \eqalign
    &\mathbb {P}_{\sinfty}^{\smallsup{\rm brw}}(\exists i_1, i_2\in T\text{ with }|i_1|=|i_2|=m, i_1\neq i_2,
    i_1, i_2\conn \infty)\nn\\
    &\qquad\leq \mathbb {P}_{\sinfty}^{\smallsup{\rm brw}}(\exists i_1, i_2\in T\text{ with }|i_1|=|i_2|=m, i_1\neq i_2,
    i_1, i_2\conn k)\nn\\
    &\qquad =\lim_{n\rightarrow \infty}\mathbb {Q}_{n}^{\smallsup{\rm brw}}(\exists i_1, i_2\in T\text{ with }|i_1|=|i_2|=m, i_1\neq i_2,
    i_1, i_2\conn k),
    \enalign
where we use Theorem \ref{thm-BRWeq} for the last equality.
We now continue to compute
    \eqalign
    &\mathbb {P}_{\sinfty}^{\smallsup{\rm brw}}(\exists i_1, i_2\in T\text{ with }|i_1|=|i_2|=m, i_1\neq i_2,
    i_1, i_2\conn \infty)\nn\\
    &\qquad \leq\lim_{n\rightarrow \infty}\frac{1}{\theta_n}\mathbb {P}(\exists i_1, i_2\in T\text{ with }|i_1|=|i_2|=m, i_1\neq i_2,
    i_1, i_2\conn k, \exists i_0 \text{ with } |i_0|=m, i_0\conn n)\nn\\
    &\qquad \leq \lim_{n\rightarrow \infty}\frac{1}{\theta_n}\mathbb {P}(\exists i_1, i_2\in T\text{ with }|i_1|=|i_2|=m, i_1\neq i_2,
    i_1\conn k, i_2\conn n)\nn\\
    &\qquad \leq \lim_{n\rightarrow \infty} \frac{\theta_{n-m}\theta_{k-m}}{\theta_n}=\theta_{k-m},
    \enalign
where we use the fact that the event that there exist $i_0, i_1, i_2$ with $i_1\neq i_2$ and such that
$i_0\conn n, i_1, i_2\conn k$ is contained in the event that there exists distinct $i_1, i_2$ such that
$i_1\conn n$ and $i_2\conn k$. Let $k\rightarrow \infty$ to obtain the result.
\qed

\subsection{Discussion and notes}
\label{sec-discBRW}
Our results for incipient infinite branching random walk are not new.
For example, the IIBRW measure in \refeq{BRWIICaltdef} was already
constructed by Kesten \cite{Kest86}, who also identified the
IIBRW measure in \refeq{IIBRWid}. See the notes below for more
details of Kesten's work. We now give an account of the history
of the problem.

IIBRW corresponds to simply embedding critical branching process trees
conditioned on non-extinction in $\Z^d$. Critical branching process trees
conditioned on non-extinction have a long history, and many properties
are known for such processes. For example, the immortal particle picture
is present in \cite[Page 304]{Geig99}, which constructs conditioned
branching processes simplifying constructions in \cite{Kest86, LN68, LPP95}.
The description in \refeq{Vzetalaw} follows from \cite[Lemma 2.1]{Geig99}.
The asymptotics in \refeq{qnasy} is shown in \cite[Theorem 3.1]{Geig99},
and goes in its most basic form back to Kolmogorov \cite{Kolm38}.
See also the references in \cite{Geig99} for a more detailed
account of the history of \refeq{qnasy}.

In \cite{Haas94}, critical branching process trees
conditioned on non-extinction are constructed for binomial offspring
distributions. This corresponds to the incipient infinite clusters
for percolation on a tree. We will discuss incipient infinite percolation
clusters in more detail in Section \ref{sec-iicop} and \ref{ss:pe}.

The immortal particle picture and the fact that there is a unique
infinite path for critical branching process trees conditioned on
non-extinction, are essential ingredients in \cite{Kest86},
where Kesten investigates random walk on a critical branching
process conditioned to survive forever. Since the probability
that the size of the total progeny for critical branching processes
exceeds $n$ decays as $1/\sqrt{n},$ these critical branching
processes create large dead ends. The random walker spends considerable
time to get out of these dead ends, which slows the random walk
down considerably, and, as a result, the graph distance of the
walker after $n$-steps grows subdiffusively as $n^{1/3}$. If we were to
embed the tree into $\Z^d$, this suggests that random walk
on the IIBRW has displacement of the order $n^{1/6}$.

The use of moment measures to describe BRW is not so common in the probability
literature. One reason may be that it they are harder to use for branching
laws that have all moments. For instance, the approach in Theorem \ref{prop-taur}
is restricted to measures having all moments. The convergence to SBM, as described in
the next section, also holds when, say, the third moment is finite.
It may be possible to use moment measures in combination with Laplace
transforms to overcome these problems. In the probability community, other
methods, such as martingale methods (see e.g.,
\cite{Perk02}) are used to prove convergence towards SBM. Such methods,
however, are based on the fact that BRW has independent branching and
motion. The models we will discuss in the sequel are self-interacting, and
moment measures are a more robust way to investigate these models.

We have proved that two constructions for IIBRW agree,
namely, (1) by size-biasing with respect to $N_n$
and letting $n\rightarrow \infty$, and (2)
by conditioning on $N_n>0$ and letting $n\rightarrow \infty$.
We believe that there are many more constructions leading
to the same limit. Examples are (3) Conditioning super-critical
branching random walk on non-extinction, and taking the limit
when the parameter turns to the critical value; (4) Conditioning
the tree $T$ to have total size $n$ and taking the limit
$n\rightarrow \infty$. It would be of interest to investigate
these constructions, and possibly other related constructions,
in more detail.

\section{Infinite canonical super-Brownian motion}
\label{sec-ICSBM}
In this section, we construct the incipient infinite canonical
measure for super-Brownian motion, which we will abbreviate as
\ICSBM (\icsbm\!\!). We will present two constructions, mirroring
the two constructions of IIBRW. We will also motivate the constructions
and definitions using the discussion of branching random walk in the
previous section to make the details comprehensible for non-specialists
in the field of super-processes.

\subsection{Super-Brownian motion and the canonical measure}
We first introduce some notation. We denote by $\N_0$ the
{\it canonical measure of super-Brownian motion}.
This canonical measure is a measure on continuous paths from
$[0,\infty)$ into non-negative finite measures on $\Rd$.
The canonical measure is an elusive object, as is it not a
{\it probability measure}, but rather a $\sigma$-finite, non-negative
measure.
We take $\N_0$ to be normalised to have unit branching
and diffusion rates. We will now first discuss a construction
of the canonical measure as a scaling limit of branching random walk
to explain the canonical measure in more detail. For
simplicity, we take an offspring distribution for which
$\sigma_p=1$. Let $n\geq 0$, and recall the definition of the random measures
$\{\mu_n\}_{n=0}^{\infty}$ in \refeq{mundef}--\refeq{mu0def}.
We expect that, as $n\rightarrow \infty$, the process $\{\mu_n\}_{n=0}^{\infty}$
has a scaling limit. The difficulty in describing this scaling limit, however, is that
    \eq
    \Pbold^{\smallsup{\rm brw}}(\exists x\in \Z^d \text{ such that }\mu_n(x)\neq 0)
    = \Pbold^{\smallsup{\rm brw}}(\exists i\in T: |i|=n) =\theta_n,
    \en
so that by \refeq{qnasy}, with probability close to 1 for $n$ large,
the random measure $\mu_n$ has mass zero. We are interested in the scaling limit,
and, in particular, in {\it large} realizations of $T$ for which
$\mu_n$ is not identically equal to 0. We now describe the construction of
the scaling limit in detail.

We define, for $t\geq 0$ and $x\in \R^d$, the random measure-valued Markov process
    \eq
    \lbeq{Xntdef}
    X_{n,t}(x)=\frac 1{n} \mu_{\lfloor nt\rfloor}(\lfloor x\sqrt{\sigma^2 n}\rfloor),
    \en
where $\lfloor x\sqrt{\sigma^2 n}\rfloor=
(\lfloor x_1\sqrt{\sigma^2 n}\rfloor, \ldots, \lfloor x_d\sqrt{\sigma^2 n}\rfloor)$.

We first motivate the scaling in \refeq{Xntdef}. It turns out that when
there is a particle alive at time $\lfloor nt\rfloor$, then there are in fact
{\it many} particles alive at the same time. Indeed, it can be shown that
conditionally on $N_m\geq 1$, the random variable $\frac{N_m}{m}$ weakly converges
to an exponential random variable. See \cite{Yagl47} or
\cite[Theorem II.1.1(b)]{Perk02}. We will be particularly interested in
branching processes that are alive at time proportional to $n$, so that
we should normalise the number of particles with a factor of $\frac 1n$.
This explains the factor $\frac{1}{n}$ in \refeq{Xntdef}. To explain the
scaling in the spatial coordinate, we note that when there is
a particle present at some site $z$ at time proportional to
$n$, then this particle has arrived to $z$ by a random walk
path of length proportional to $n$. Therefore, we can expect
that $z$ is of the order $\sqrt{\sigma^2 n}$. This explains
the scaling in \refeq{Xntdef}. We now describe the scaling limit.

For an event $E$ that is a measurable subset of the space of measure-valued
paths on $\R^d$, we take the limit
    \eq
    \lbeq{WCCM}
    \lim_{n\rightarrow \infty} n \Pbold^{\smallsup{\rm brw}}
    (\{X_{n,t}\}_{t\geq 0} \in E).
    \en
It turns out that the above limit exists as an element of $[0,\infty]$
(see e.g., \cite[Theorem II.7.3(a)]{Perk02}), and is
by definition equal to the measure of the indicator of the
event $E$ under the canonical measure of super-Brownian motion, i.e.,
to $\N_0(I[E])$.
The factor $n$ in \refeq{WCCM} explains that the measure $\N_0$
is not a probability measure, but rather a $\sigma$-finite measure.
We now discuss this construction of the canonical measure and its
relation to super-Brownian motion started from a proper initial
measure.

For a measure $\mu$ on $\R^d$, we write $\mu(1)=\int_{\R^d} 1 d\mu$ for
its total mass. We write $\{X_t\}_{t\geq 0}$ for the process of non-negative
measures under the canonical measure $\N_0$. Note that when
$E=\{X_t(1)>0\}$, then by \refeq{qnasy},
    \eq
    \lbeq{N0t}
    \N_0(I[X_t(1)>0])=\lim_{n\rightarrow \infty} n \Pbold^{\smallsup{\rm brw}}(X_{n,t}(1)>0)
    =\lim_{n\rightarrow \infty} n \theta_{\lfloor nt\rfloor}=\frac{2}{\sigma_p^2t}
    =\frac{2}{t},
    \en
since we have assumed that $\sigma_p=1$. Therefore, $\N_0$ is a {\it finite}
measure on events $E$ that imply that $X_t(1)>0$ for some $t>0$.

Often, super-Brownian motion is considered as starting from a {\it proper
initial measure}. This corresponds to a different scaling limit. Indeed, let
the measure $\mu_{0,n}$ be such that $\mu_{0,n}(x)$ takes integer
value for every $x\in \Z^d$, and let
    \eq
    \nu_n(x) = \frac 1{n} \mu_{0,n}(\lfloor x\sqrt{\sigma^2 n}\rfloor).
    \en
We assume that $\nu_n$ is a measure that weakly converges to some
limiting measure $\nu$. Then, we let $(T_{x}^{\smallsup{j}},
\phi_{x}^{\smallsup{j}})$ for $j=1, \ldots, \mu_{0,n}(x)$
be $\mu_{0,n}(x)$ independent branching random walks
started at $x$, so that $\phi_{x}^{\smallsup{j}}(0)=x$,
where $0$ is the root of the tree $T_{x}^{\smallsup{j}}$.
We now start with several independent branching random walks
with starting points given by the initial measure $\mu_{0,n}$.
The spatial locations of the branching random walk
particles at time $m$ with initial measure $\mu_n$ are then given by
    \eq
    \mu_{m,n}(x)=\sum_{y\in \Z^d} \sum_{j=1}^{\mu_{0,n}(y)}
    \sum_{i\in T_{y}^{\smallsup{j}}: |i|=m}
    I[\phi_{y}^{\smallsup{j}}(i)=x].
    \en
In words, the random variable $\mu_{m,n}(x)$ equals the number of particles that
are present at time $m$ at the location $x$ when we start with initial
measure $\mu_{0,n}$. Then, we define
    \eq
    \lbeq{XntdefSBM}
    X_{n,t}(x)=\frac 1{n} \mu_{\lfloor nt\rfloor,n}(\lfloor x\sqrt{\sigma^2 n}\rfloor)
    \en
Thus, in particular, $X_{n,0}(x)=\nu_n(x)$.
Denote the law of $\{X_{n,t}\}_{t=0}^{\infty}$ by $\Pbold_{\nu_{n}}^{\smallsup{\rm brw}}$.
Then, the limit
    \eq
    \lbeq{WCSBM}
    \lim_{n\rightarrow \infty} \Pbold_{\nu_{n}}^{\smallsup{\rm brw}}(\{X_{n,t}\}_{t\geq 0} \in E)
    \en
exists. This limit is
    \eq
    \Pbold_{\nu}(\{X_{t}\}_{t\geq 0} \in E),
    \en
where $\Pbold_{\nu}$ is the law of {\it super-Brownian motion with
initial measure} $\nu$.

We can think of the law $\Pbold_{\nu_{n}}^{\smallsup{\rm brw}}$ as being
described by the evolution of independent branching random walk
copies, where the copies are located at positions
described by the initial measure $\nu_{n}$. In a similar way, we can
think of $\Pbold_{\nu}$ as being described by (infinitely) many
independent copies of canonical measures according to the initial
measure $\nu$ (see e.g., \cite[Theorem II.7.2]{Perk02}).
This intuitive picture can be made precise by noting that
$\Pbold_{\nu}$ is {\it infinitely divisible}, and using
the general notion of infinitely divisible measures in e.g.,
\cite{Kall83}. In fact, in the terminology of infinite divisible measures,
the canonical measure of super-Brownian motion is the canonical
measure for the infinitely divisible measure $\Pbold_{\nu}$.
See e.g.\ \cite[Section 1.3, and, in particular, Corollary 1.3]{DaP99}
for more details.

On the other hand, we can also describe the canonical measure
in terms of SBM by using the Markov property.
Indeed, the law of $\{X_{s+t}\}_{t=0}^{\infty}$
given $X_s$ is the same as the law of $\{X_{t}\}_{t=0}^{\infty}$
under $\Pbold_{X_s}$.\ \footnote{This follows immediately from the Markov
property valid for BRW, together with the weak convergence in \refeq{WCCM}
and \refeq{WCSBM}. See also \cite[Theorem 1.4 (iii)]{DaP99}, where the Markov
property is stated for the canonical measure.}
This shows that the laws of SBM with a proper
initial measure and the canonical measure are
intimately connected.


\subsection{The infinite canonical measure of super-Brownian motion}
We are now ready for the definition of \icsbm\!\!, which we will think of as
the canonical measure conditioned on non-extinction.
Fix $t>0$. Let $\Fcal_s$ be the minimal $\sigma$-algebra
of events such that $\{X_u\}_{0\leq u\leq s}$ is measurable with respect to
$\Fcal_s$, i.e., $\Fcal_s=\sigma(\{X_u\}_{0\leq u\leq s})$, and we write
    \eq
    \Fcal=\bigcup_{s\geq 0} \Fcal_s.
    \en
Let $E\in\Fcal_s$ for some $s<\infty$. Thus,
$E$ only depends on $X_u$ for all $u\leq s$. We then define the
probability measure
    \eq
    \lbeq{Ptdef}
    \Pbold_t(E) = \frac{\N_0(I[E] X_t(1))}{\N_0(X_t(1))},
    \en
where the random variable $X_t(1)$ is the total mass of super-Brownian motion
at time $t$, and where, for an event $E$, $I[E]$ is the indicator of the event
$E$. The measure $\Pbold_t$ is {\it size-biased} with respect to
the total mass at time $t$. The definition in \refeq{Ptdef} is reminiscent
of the definition in \refeq{PnBRWdef}.

We also define a second probability measure
    \eq
    \lbeq{Qtdef}
    \Qbold_t(E) = \frac{\N_0(I[E] I[X_t(1)>0])}{\N_0(I[X_t(1)>0])},
    \en
which we can think of as the canonical measure of
super-Brownian motion conditioned to survive up to
time $t$, and which is reminiscent of the definition
in \refeq{QnBRWdef}.

The measures $\Pbold_t$ and $\Qbold_t$, for large $t$, ensure that
the SBM does not have too small mass. Indeed, when $t$ is large, then
with probability close to 1, we have that $X_t(1)=0$. On the contrary, for large $t$,
SBM must live until time $t$ in both constructions. We now prove that
the above two measures converge, when $t\rightarrow \infty,$ to a
limiting measure, which we will call {\it \ICSBM\!\!:}

    \begin{theorem}
    \label{thm-ISBM}
    When $t\rightarrow \infty$, and for every $E\in \Fcal_s$ for any $s\geq 0$,
    $\Pbold_t(E)$ and $\Qbold_t(E)$ converge to
    $\Pbold_{\sinfty}(E)$ and $\Qbold_{\sinfty}(E)$.
    Moreover, $\Pbold_{\sinfty}$ and $\Qbold_{\sinfty}$ extend to
    probability measures on the $\sigma$-algebra
    $\Fcal$, and $\Pbold_{\sinfty}=\Qbold_{\sinfty}.$
    \end{theorem}

\proof
The statement in Theorem \ref{thm-ISBM} is simplest to prove for $\Pbold_{\sinfty}$. Indeed,
since the total mass $X_t(1)$ is a martingale with $X_0(1)=1$, we have that
(see also \cite[Theorem II.7.2 (iii)]{Perk02})
    \eq
    \N_0(X_t(1))=1,
    \en
and, when $E\in\Fcal_s$,
    \eq
    \N_0(I[E] X_t(1))=\N_0(I[E] X_s(1)).
    \en
Thus, when $E\in\Fcal_s$,
    \eq
    \lbeq{Ptseq}
    \Pbold_t(E) = \Pbold_s(E),
    \en
and the stated convergence is trivial. Since the limiting measure $\Pbold_{\sinfty}$
is consistent, we can extend it to the full $\sigma$-algebra by Kolmogorov's
Extension Theorem (see e.g.\ \cite{Simo79}). For the second construction, we have to
do a bit more work. We note from \cite[Theorem II.7.2 (iii)]{Perk02} that
    \eq
    \N_0(I[X_t(1)>0]) = \frac{2}{t},
    \en
(see also \refeq{N0t} above) so that it suffices to compute
$\N_0(I[E] I[X_t(1)>0])$. We will first prove that
    \eq
    \lbeq{convstat}
    \lim_{t\rightarrow \infty} t\mathbb{P}_{\mu} (X_t(1)>0)=2\mu(1),
    \en
where $\mathbb{P}_{\mu}$ is the law of super-Brownian motion starting from the
measure $\mu$. Equation \refeq{convstat} follows from the fact that
    \eq
    \lbeq{convstat2}
    \mathbb{P}_{\mu} (X_t(1)>0) = 1-e^{-2\mu(1)t^{-1}}.
    \en
(See \cite[(1.3)]{EP90}, where the factor 2 is absent due to the fact we have assumed the
branching rate to be 1, whereas there it is 2. See also \cite[(II.5.12)]{Perk02}.)
Using \refeq{convstat}, we can write
    \eq
    t\N_0(I[E] I[X_t(1)>0]) = \N_0(I[E] t\mathbb{P}_{X_s}(X_{t-s}(1)>0)).
    \en
Here we use the {\it Markov property for the canonical measure} described above.
Using \refeq{convstat2}, and the fact that thus, for $t\geq 2s$,
    \eq
    t\mathbb{E}_{X_s}[X_{t-s}(1)>0]\leq 2X_s(1)\frac{t}{t-s} \leq 4X_s(1),
    \en
which is integrable, we obtain by dominated convergence,
    \eq
    \lbeq{Qtcomp}
    \lim_{t\rightarrow \infty}\Qbold_t(E) =
    \lim_{t\rightarrow \infty}\frac{\N_0(I[E] I[X_t(1)>0])}{\N_0(I[X_t(1)>0])}
    =\lim_{t\rightarrow \infty}\frac{t}{2}\N_0(I[E] I[X_t(1)>0]) = \N_0(I[E]X_s(1))
    =\mathbb{P}_{\sinfty}(E).
    \en
This completes the proof of Theorem \ref{thm-ISBM}.
\qed

\subsection{Convergence of BRW $r$-point functions to SBM moment measures}
In this section, we describe the moment measures of the canonical measure of
super-Brownian motion. We first discuss what these moment measures are and how they can
be characterized. In the next section, we use this characterization to describe the
\icsbm moment measures.

As we will explain in more detail below, moment
measures describe the finite-dimensional distributions of super-process.
Indeed, a measure can be determined by its expectation of a sufficiently
rich class of bounded continuous functions. For a random measure $X_t$, we can
thus determine the law of $X_t$ by describing the laws of $X_t(f)$ for a
sufficiently rich class of continuous functions, where
    \eq
    X_t(f) = \int_{\R^d} f(x) dX_t(x).
    \en
We will be using Fourier transforms, so that we take as a class of
continuous functions $\{f_k\}_{k \in \R^d}$, where $f_k(x)=e^{ik\cdot x},$
and $k\cdot x$ is the inner product between $x$ and $k$. Thus, in order to determine the law of
super-Brownian motion, it suffices to know the law of $\{X_s(f_k)\}_{s\geq 0, k\in \R^d}$.
This law will be uniquely determined by the finite-dimensional distributions
$\{X_{s_i}(f_{k_i})\}_{i=1}^r$ for any $s=(s_1, \ldots, s_r)$ and $k=(k_1, \ldots, k_r)$.
These laws, in turn, will be unique determined in terms of the joint moments,
for every vector $(a_1, \ldots, a_r)\in \N^r$,
    \eq
    \N_0\big(\prod_{i=1}^r X_{s_i}(f_{k_i})^{a_i}\big)=\N_0\big(\prod_{j=1}^l X_{t_j}(f_{k_j})\big),
    \en
where $l=a_1+\ldots+a_r$, and the components of $(t_1, \ldots, t_l)$ are equal
to $s_j$ precisely $a_j$ times. Thus, we are lead to investigate
    \eq
    \hat{M}^{\smallsup{l}}_{\vec{t}}(\vec{k})=\N_0\big(\prod_{j=1}^l X_{t_j}(f_{k_j})\big)=
    \N_0\big(\int_{\R^{dl}}X_{t_1}(dx_1) \cdots X_{t_{l}}(dx_{l})
    \prod_{j=1}^{l} e^{ik_j\cdot x_j}\big).
    \en
These are the Fourier transforms of the
{\it moment measures of the canonical measure of
super-Brownian motion}.

We next give formulas for these moment measures of the
canonical measure of super-Brownian motion. We follow
the presentation in \cite[Section 4]{HHS02a}. We will make use of
elementary properties of the $\hat{M}^{\smallsup{l}}_{\vec{t}}(\vec{k})$,
which we now summarise.
For $l=1$,
    \eq
    \lbeq{fdSBM}
    \hat{M}^{\smallsup{1}}_{t}(\vec{k})
    =
    e^{-|k|^2t/2d}.
    \en
We will write
    \eq
    \lbeq{knotation}
    k_I=\sum_{i\in I} k_i, \qquad \vec{k}_{I}= (k_i; i\in I).
    \en
For $l>1$, the $\hat{M}^{\smallsup{l}}_{\vec{t}}(\vec{k})$
are given recursively by
    \eq
    \lbeq{recmeas}
        \hat{M}_{\vec{t}}^{\smallsup{l}}(\vec{k})
        =
        \int_0^{\underline{t}} dt \;
        \hat{M}_{t}^{\smallsup{1}}(k_J)
        \sum_{I \subset J_1 : |I| \geq 1}
        \hat{M}_{\vec{t}_{I}-t}^{\smallsup{i}}(\vec{k}_{I})
        \hat{M}_{\vec{t}_{J\backslash I}-t}^{\smallsup{l-i}}
        (\vec{k}_{J\backslash
        I}),
    \en
where $i=|I|$, $J=\{1, \ldots, l\}, J_1=J\backslash \{1\}$,
$\underline{t}=\min_i t_i$, $\vec{t}_I$ denotes the vector consisting
of the components $t_i$ of $\vec{t}$ with $i \in I$, and
$\vec{t}_{I}-t$ denotes subtraction of $t$ from each component
of $\vec{t}_I$ \cite{Dynk88}.  The explicit solution to the recursive formula
\refeq{recmeas} can be found in \cite[(1.25)]{HS02a}.  For example,
    \eq
    \lbeq{M2}
    \hat{M}_{t_1,t_2}^{\smallsup{2}}(k_1,k_2)
    =
    \int_0^{t_1 \wedge t_2}
    dt \;
    e^{-|k_1+k_2|^2t/2d} e^{-|k_1|^2(t_1-t)/2d} e^{-|k_2|^2(t_2-t)/2d}.
    \en
Equation~\refeq{M2} is a statement, in Fourier language, that mass arrives
at given points $(x_1,t_1)$, $(x_2,t_2)$ via a Brownian path from
the origin that splits into two Brownian paths at a time chosen uniformly
from the interval $[0,t_1 \wedge t_2]$.  The recursive formula
\refeq{recmeas} has a related interpretation for all $l \geq 2$,
in which $t$ is the time of the first branching.  The sets $I$ and $J
\backslash I$ label the offspring of each of the two particles after
the first branching.

The main result in this section is the proof that the moment measures of
the canonical measure of SBM arise as the scaling limits of the BRW
$r$-point functions:
\begin{theorem}
\label{prop-condBRW}
Fix an offspring distribution $(p_m)_{m=0}^{\infty}$ such that
all moments are finite, and assume that \refeq{deltadef} holds. Then,
   \eq
   \hat{\tau}_{\vec{n} }(\vec{k}/\sqrt{\sigma^2 n})
    =(\sigma_p^{2}n)^{r-2}
    \left[
    \hat{M}_{\vec{n}/n}^{\smallsup{r-1}}(\vec{k})
    + {\cal O}((n_{\smallsup{2}}+1)^{-\delta})
        \right]
    \quad
    (r \geq 2)
    \lbeq{convSBMBRW}
    \en
holds uniformly in $n\geq n_{\smallsup{2}}$, where $n_{\smallsup{2}}$ denotes the
second largest component of $\vec{n}$.
In particular, $\sigma_p^2 n\Pbold^{\smallsup{\rm brw}}(X_{n,\frac{\vec{t}}{\sigma_p^2}} \in \cdot)$ converges to
$\N_0$ in the sense of convergence of finite-dimensional distributions.
\end{theorem}

\proof
We start from the recursive formula in Proposition \ref{prop-taur}, and rewrite it as
    \eq
    \lbeq{taurew}
    \tau_{\vec n}(\vec x)=(D*\tau_{\vec n-1})(\vec{x})
    +\sigma_p^2  \sum_{I\not\ni 1, I\neq \varnothing}
    (D*\tau_{\vec n_{I}-1})(\vec{x}_{I})
    (D*\tau_{\vec n_{J\backslash I}-1})(\vec{x}_{J\backslash I})
    +e_{\vec n}(\vec x),
    \en
where
    \eq
    e_{\vec n}(\vec x)=\sum_{j=3}^{r-1} f_j \sum_{\vec{I}\in \Pcal_j}
    \prod_{s=1}^{r-1} (D*\tau_{\vec n_{I_s}-1})(\vec{x}_{I_s}).
    \en
Thus, $e_{\vec n}(\vec x)$ is the contribution to $\tau_{\vec n}(\vec x)$ where there are
at least three children of the root that are connected to elements $i_j$ with $|i_j|=n_j$ and $\phi(i_j)=x_j$.
Equation \refeq{taurew} is a generalization of \refeq{taureq}, which was only valid for
binary branching. In deriving \refeq{taurew}, we have use that $f_1=1, f_2=\sigma_p^2$.

We will prove that when $(p_m)_{m=0}^{\infty}$ has all moments, then there exists
$C_r<\infty$ such that
    \eq
    \lbeq{enbd}
    |\hat{e}_{\vec n}(\vec k)|\leq \sum_{\vec{x}} e_{\vec n}(\vec x)\leq C_r (n_{\smallsup{2}}+1)^{r-3}.
    \en
We first show that this suffices to prove the statement in Proposition \ref{prop-condBRW}.
Indeed, we iterate \refeq{taurew} until the first term disappears. Then we arrive at
    \eq
    \lbeq{taurew2}
    \tau_{\vec n}(\vec x)=\sigma_p^2  \sum_{m=0}^{\underline{n}}\sum_{I\not\ni 1, I\neq \varnothing}
    \sum_y
    D^{*m}(y)
    (D*\tau_{\vec n_{I}-m-1})(\vec{x}_{I}-y)
    (D*\tau_{\vec n_{J\backslash I}-m-1})(\vec{x}_{J\backslash I}-y)
    +\varphi_{\vec n}(\vec x),
    \en
where
    \eq
    \varphi_{\vec n}(\vec x) =\sum_{m=0}^{\underline{n}} \sum_y
    D^{*m}(y) e_{\vec n-m}(\vec x-y).
    \en
Clearly, by \refeq{enbd},
    \eq
    |\hat{\varphi}_{\vec n}(\vec k)|\leq \sum_{\vec{x}} \varphi_{\vec n}(\vec x)\leq C_r n_{\smallsup{2}}^{r-2},
    \en
so that $|\hat{\varphi}_{\vec n}(\vec k)|$ is an error term. We take the Fourier
transform of \refeq{taurew2} to obtain
    \eq
    \lbeq{taurew3}
    \hat{\tau}_{\vec n}(\vec k)
    =
    \sigma_p^2  \sum_{m=0}^{\underline{n}}\sum_{I\not\ni 1, I\neq \varnothing}
    \hat D^{m}(k_J)
    \hat D(k_I)\hat\tau_{\vec n_{I}-1}(\vec{k}_{I})
    \hat D(k_{J\backslash I}) \hat\tau_{\vec n_{J\backslash I}-1}(\vec{k}_{J\backslash I})
    +\hat\varphi_{\vec n}(\vec k).
    \en
Equation \refeq{taurew3} is a discrete version of \refeq{recmeas}, and it can be used
to prove by induction on $r$ that
    \eq
    \hat{\tau}_{\vec{n} }(\vec{k}/\sqrt{\sigma^2 n})
    =  \sigma_p^{2(r-2)} n^{r-2}
    \left[
    \hat{M}_{\vec{n}/n}^{\smallsup{r-1}}(\vec{k})
    + {\cal O}((n_{\smallsup{2}}+1)^{-\delta})
        \right]
    \quad
    (r \geq 2)
    \lbeq{aim3ptBRW}
    \en
holds uniformly in $n\geq n_{\smallsup{2}}.$ For $r=2$, we use \refeq{tauneq}, which implies
that $\hat{\tau}_n(k/\sqrt{\sigma^2 n})=\hat{D}(k/\sqrt{\sigma^2 n})^n$.
When \refeq{deltadef} holds, then $\hat{D}(k)=1-\sigma^2 \frac{|k|^2}{2d} +\Ocal(|k|^{2+2\delta}$.
This immediately implies the claim for $r=2$, and initializes the induction.
We will omit the details of the advancement of the induction, which can be found
in \cite[Section 2.3]{HS02a}, where the same computation was performed for oriented percolation.
This completes the proof of Proposition \ref{prop-condBRW} subject to \refeq{enbd}. We complete the
proof by proving \refeq{enbd}.

We use \refeq{taurew2} to prove by induction on $r$ that there exists a $C_r<\infty$
such that
    \eq
    \hat \tau_{\vec{n}}(\vec{0})=\sum_{\vec x} \tau_{\vec{n}}(\vec{x})\leq C_r (n_{\smallsup{2}}+1)^{r-2}.
    \en
For $r=2$, the inequality holds, since the left-hand side equals 1. This initialises
the induction.

Without loss of generality, we may assume that $C_r$ is non-decreasing in $r$,
and that $C_r\geq 1$. To advance the induction, we note that when the above claim holds
for all $s\leq r-1$, then, using \refeq{taurew2}, we obtain that
    \eqalign
    \hat \tau_{\vec{n}}(\vec{0}) &\leq  \sum_{j=2}^{r-1}
    \sum_{m=0}^{\underline{n}}
    f_j \sum_{\vec{I}\in \Pcal_j}
    \prod_{s=1}^{r-1} \hat \tau_{\vec n_{I_s}-1}(\vec{0}_{I_s})
    \leq  C_{r-1}^{r-1}\sum_{j=2}^{r-1}
    f_j \sum_{\vec{I}\in \Pcal_j}
    \sum_{m=0}^{\underline{n}} \prod_{s=1}^{r-1} (n_{\smallsup{2}}-m)^{|I_s|-1}\nn\\
    & \leq C_{r-1}^{r-1}\sum_{j=2}^{r-1}
    f_j \sum_{\vec{I}\in \Pcal_j} (n_{\smallsup{2}}+1)^{\sum_{s=1}^{r-1} |I_s|-j+1}
    \leq C_r (n_{\smallsup{2}}+1)^{r-2},
    \enalign
where we use that $|I_1|+\ldots+|I_{j}|=r-1$, and where $C_r$
must be chosen appropriately large. A similar computation, where the sum over $j$
starts at $j=3$, proves \refeq{enbd}.
\qed

\subsection{The moment measures of \icsbm}
\label{sec-MMSBM}
In this section, we describe the moment measure of \icsbm\!\!.

We write $M_{\sinfty;\vec{s}}^{\smallsup{l}}$ for the moment measures of \icsbm\!\!,
and $\hat{M}_{\sinfty;\vec{s}}^{\smallsup{l}}(\vec{k})$ for the Fourier
transform of the \icsbm moment measures, i.e.,
\eq
    \hat{M}^{\smallsup{l}}_{\sinfty;\vec{s}}(\vec{k})
        =\Ebold_{\sinfty}\big(\int_{\R^{dl}} X_{s_1}(dx_1) \cdots X_{s_{l}}(dx_{l})
    \prod_{j=1}^{l} e^{ik_j\cdot x_j}\big),
\en
where $\vec{s}=(s_1,\ldots,s_l)$ with each $s_i\in (0,\infty)$, and
$\vec{k}=(k_1,\ldots,k_l)$ with each $k_i \in \Rd$.

We can then identify the moment measures of \icsbm as follows:

\begin{theorem}
\label{thm-MMISBM}
For every $\vec{s}=(s_1,\ldots,s_l)$ with each $s_i\in (0,\infty)$, and
$\vec{k}=(k_1,\ldots,k_l)$ with each $k_i \in \Rd$,
    \eq
    \lbeq{MPSBM}
    \hat{M}_{\sinfty;\vec{s}}^{\smallsup{l}}(\vec{k})=
    \hat{M}_{\bar{s}, \vec{s}}^{\smallsup{l+1}}(0, k_1, \ldots, k_l),
    \en
where $\bar{s}=\max_{1\leq i\leq l} s_i$.
   \end{theorem}

\proof
A similar statement as in \refeq{MPSBM} was proved in \cite[Lemma 4.2]{HHS02a}
when $s_i=s$ for some $s$. This proof uses induction on $l$. We now prove
the more general version of this claim using a simpler martingale proof.
We note that the integral
    \[\int_{\R^{dl}} X_{s_1}(dx_1) \cdots X_{s_{l}}(dx_{l})
    \prod_{j=1}^{l} e^{ik_j\cdot x_j}
    \]
only depends on $X_s$ for $s\leq \bar{s}$.
Therefore, by Theorem \ref{thm-ISBM},
    \eq
    \hat{M}^{\smallsup{l}}_{\sinfty;\vec{s}}(\vec{k})
        =\lim_{t\rightarrow \infty}\Ebold_{t}\big(\int_{\R^{dl}} X_{s_1}(dx_1) \cdots X_{s_{l}}(dx_{l})
    \prod_{j=1}^{l} e^{ik_j\cdot x_j}\big).
    \en
By \refeq{Ptseq}, we have that
    \eqalign
    \Ebold_{t}\big(\int_{\R^{dl}} X_{s_1}(dx_1) \cdots X_{s_{l}}(dx_{l})
    \prod_{j=1}^{l} e^{ik_j\cdot x_j}\big)
    &=\Ebold_{\bar{s}}\big(\int_{\R^{dl}} X_{s_1}(dx_1) \cdots X_{s_{l}}(dx_{l})
    \prod_{j=1}^{l} e^{ik_j\cdot x_j}\big)\nn\\
    &=\Ebold\big(\int_{\R^{dl}} X_{s_1}(dx_1) \cdots X_{s_{l}}(dx_{l})X_{\bar{s}}(1)
    \prod_{j=1}^{l} e^{ik_j\cdot x_j}\big)\nn\\
    &=\hat{M}_{\bar{s}, \vec{s}}^{\smallsup{l+1}}(0, k_1, \ldots, k_l).
    \enalign
for every $t\geq \bar{s}$. This completes the identification of the
moment measures of \icsbm\!\!.
    \qed

Note that, in particular,
    \eq
    \lbeq{MM1}
    \hat{M}^{\smallsup{1}}_{\sinfty;s}(k)
        =\hat{M}_{s,s}^{\smallsup{2}}(0,k)=s e^{-\frac{k^2s}{2d}},
    \en
and
    \eq
    \lbeq{MM2}
    \hat{M}^{\smallsup{2}}_{\sinfty;\vec s}(\vec k)
        =\int_0^{s_1 \wedge s_2} (s_1+s_2-s)
        e^{-|k_1+k_2|^2s/2d} e^{-|k_1|^2(s_1-s)/2d} e^{-|k_2|^2(s_2-s)/2d}ds.
    \en

We next investigate the total mass under the ICSBM measure.
\begin{theorem}
\label{thm-sbExp}
For $s>0$, $X_s(1)$ is
a size-biased exponential random variable with
parameter $2/s.$
\end{theorem}

\proof In \cite[Lemma 4.2(c)]{HHS02a}, it was proved that for $l \geq 0$,
    \eq
    \lbeq{Mform}
    \Ebold_{\sinfty} [X_s(1)^l]=\hat{M}_{s,\ldots, s}^{\smallsup{l+1}}(\vec{0})= \N_0[X_s(1)^{l+1}]= s^l 2^{-l} (l+1)!,
    \en
which are the moments of a size-biased exponential random variable.
The distribution of the size-biased exponential random variable
is determined by its moments, since its moment generating function
has a positive radius of convergence (see \cite[Theorem~30.2]{Bill95}).
It therefore follows from the moments in \refeq{Mform} that $X_s(1)$ is
a size-biased exponential random variable with
parameter $2.$
\qed

We note that the exact equality in law is due to the fact
that we start with the canonical measure of super-Brownian motion. Indeed,
in \cite[Theorem (iii)]{EP90}, it follows that for a general measure-valued
process with a starting measure $\mu$, and conditionally on survival at
time $s$, the random variable $X_s(1)/s$ converges weakly
to a size-biased random variable, rather than being precisely equal to it for all
$s$.

\subsection{Convergence of the IIBRW moment measures}
We now turn to the scaling limit of incipient infinite
branching random walk. Given the close connections between SBM and
critical BRW, it can be expected that the $r$-point functions for IIBRW
converge, appropriately scaled, to their continuous analogues for SBM.
This is not trivial, since it involves the interchange of the limits defining
IIBRW and \icsbm and the scaling limit for IIBRW. The main result is the
following theorem:

    \begin{theorem}
    \label{thm-IBRW}
    \lbeq{ISBMrpt}
    Fix an offspring distribution $(p_m)_{m=0}^{\infty}$ such that
    all moments are finite. Then, for all $r \geq 2$, $\vec{t} =
    (t_1,\ldots,t_{r-1}) \in \mathbb{R}_+^{r-1}$
    and $\vec{k} \in \R^{d(r-1)}$, and with $\delta\in (0,1)$ as
    in \refeq{deltadef},
    \eq
    \lbeq{rhoscal}
    \frac{1}{(A^2V)^{r-1}}
    \hat{\rho}_{m\vec{t}}(\vec{k}/\sqrt{\sigma^2 m})
    =
    m^{r-1}\hat{M}_{\sinfty,\vec{t}}^{\smallsup{r}}(\vec{k})[1+\Ocal(m^{-\delta})].
    \en
    Consequently, when time is rescaled by $m$ and space by $\sqrt{\sigma^2 m}$,
    the finite-dimensional distributions of ${\Pbold}_{\sinfty}^{\smallsup{\rm brw}}$
    converge to those of ${\Pbold}_{\sinfty}$.
    \end{theorem}

  \proof We prove a more general version of this result, which we can apply
  later on to other incipient infinite structures. We will prove the following proposition.
  We will assume that $\tau$ and $\rho$ are functions that are related via (recall \refeq{relrhotaubrw})
   \eq
    \lbeq{relrhotau}
    \rho_{n_1,\ldots,n_{r-1}}(x_1,\ldots, x_{r-1})
    =\lim_{n\rightarrow \infty} \sum_{x_0\in \Z^d}\tau_{n, n_1,\ldots,n_{r-1}}
    (x_0,x_1,\ldots, x_{r-1}).
    \en

  \begin{prop}
  \label{prop-cond}If there exist constants $A,V, v, \delta$ such that
   \eq
   \hat{\tau}_{\vec{n} }(\vec{k}/\sqrt{v\sigma^2 n})
    = A (A^2V)^{r-2} n^{r-2}
    \left[
    \hat{M}_{\vec{n}/n}^{\smallsup{r-1}}(\vec{k})
    + {\cal O}((n_{\smallsup{2}}+1)^{-\delta})
        \right]
    \quad
    (r \geq 2)
    \lbeq{aim3pt}
    \en
    holds uniformly in $n\geq n_{\smallsup{2}}$, then
    \eq
    \frac{1}{(A^2V)^{r-1}}
    \hat{\rho}_{m\vec{t}} (\vec{k}/\sqrt{v\sigma^2 m})
    =
    m^{r-1}\hat{M}_{\sinfty,\vec{t}}^{\smallsup{r}}(0,\vec{k})[1+\Ocal(m^{-\delta})].
    \lbeq{rhoscal2}
    \en
    \end{prop}

    \noindent
    {\it Proof of Proposition \ref{prop-cond}.}
    Equation \refeq{rhoscal} is an immediate consequence of the assumption \refeq{aim3pt}, together with the
    relation in \refeq{relrhotau}.
    \qed

The proof of Theorem \ref{thm-IBRW} is now a combination of Proposition \ref{prop-cond} and
Theorem \ref{prop-condBRW}.

We close this section with a description of the number of particles alive
at time $m$ as a corollary to Theorem \ref{thm-IBRW}.
For this, we recall that the {\em size-biased}\/ exponential
random variable with parameter $\lambda$ has density
    \eq
    f(x) = \lambda^2 x e^{-\lambda x} \qquad (x\geq 0).
    \en

\begin{theorem}
\label{thm-NmBRW}
Under $\mathbb {P}_{\sinfty}^{\smallsup{\rm brw}}$, $\frac{N_m}{m}$ converges weakly
to a size-biased exponential random
variable with parameter $\lambda=\frac{2}{\sigma_p^2}$.
\end{theorem}

\proof We compute the moments of $\frac{N_m}{m}$ under the measure
$\mathbb {P}_{\sinfty}^{\smallsup{\rm brw}}$, which equal
    \eq
    m^{-l}\Ebold_{\sinfty}^{\smallsup{\rm brw}}[N_m^l]
    =m^{-l}\hat{\rho}_{m, \ldots, m}(\vec{0}).
    \en
By Theorem \ref{thm-IBRW}, the right-hand side converges to $\sigma_p^{2l}\hat{M}_{\sinfty,\vec{1}}^{\smallsup{l}}(\vec{0})$,
which is the $l^{\rm th}$-moment of a size-biased exponential random variable (see the proof of
Theorem \ref{thm-sbExp}). Then we can follow the remainder of the proof of Theorem \ref{thm-sbExp}.
\qed

\subsection{The four-dimensional nature of \icsbm}
We next describe the four-dimensional nature of \icsbm\!\!.
Let
    \eq
    \lbeq{MRdef}
    M(R) = \int_{0}^{\infty} X_s(B_R) ds,
    \en
where $B_R$ is the unit ball of radius $R$. Then  we have the following scaling result:

    \begin{theorem}
    \label{thm-4d}Under $\Pb_{\sinfty}$, $M(R)$ has the same law as $R^4M(1)$.
    Moreover, when $d>4$, $M(1)<\infty$ $\Pb_{\sinfty}$-a.s.
    \end{theorem}

   \proof We use the Brownian scaling, which implies that
   $\{R^{-2}X_{R^2t}(\cdot R)\}_{0\leq s<\infty}$ has the same law as
    $\{X_{t}(\cdot)\}_{0\leq s<\infty}$. This equality in law can for instance
    be seen by comparing the moment measures of $\{R^{-2}X_{R^2t}(\cdot R)\}_{0\leq s<\infty}$,
    and by proving that these are equal to the ones of $\{X_{t}(\cdot)\}_{0\leq s<\infty}$.

    Therefore,
    \eqalign
    M(R) &= \int_{0}^{\infty} X_s(RB_1) ds
    =R^2\int_{0}^{\infty} X_{R^2s}(RB_1)ds
    \stackrel{\sss d}{=}R^2\int_{0}^{\infty} R^2 X_{s}(B_1)ds=R^4 M(1).
    \enalign
    The fact that $M(1)$ is finite a.s.\ when $d>4$ can be deduced from the fact that
    $\Ebold_{\sinfty}[M(1)]<\infty$, which follows from the fact that
        \eq
        \Ebold_{\sinfty}[M(1)] = \Ebold_{\sinfty} \big[\int_{0}^{\infty} X_{s}(B_1)ds\big]
        =\int_{0}^{\infty} \Ebold_{\sinfty}[X_{s}(B_1)]ds
        =\int_{0}^{\infty} \int_{B_1} s p_s(y)dy ds <\infty,
        \en
    where we used that
        \eq
        \Ebold_{\sinfty}[X_{s}(B_1)]=\int_{B_1} s p_s(y)dy,
        \en
    which follows from \refeq{MM1}
    and the Fourier-inversion formula.
    \qed

Theorem \ref{thm-4d} is important, as it suggests that critical structures
can only converge to \icsbm when they are four-dimensional.
It is possible to describe the law of $M(1)$ by computing its moments
using Theorem \ref{thm-MMISBM}, but we will not do so here.

\subsection{The immortal particle}
In \cite{Evan93}, among other things, there is a nice description of a super-processes conditioned
on non-extinction, which goes under the name of {\it Evan's immortal particle}. Roughly speaking, when we condition
a super-process starting from a proper starting distribution to survive up to some time
$t$, then we see only the descendants of a {\it finite} number of particles.
When we then take the limit of $t\rightarrow \infty$, this number of particles becomes
one, and we end up with the descendants of the immortal particle. Therefore, we can think of
\icsbm as being built up from a single particle undergoing
Brownian motion, which gives rise to offspring at a fixed rate.
This construction is similar in spirit to the one for IIBRW in Theorem \ref{thm-immparti},
even though the relation to BRW is not mentioned in \cite{Evan93}.

The construction in \cite{Evan93} is for a super-process starting in a proper
initial measure. Since we work with the canonical measure of super-Brownian motion,
we look at the super-process started with a {\it single} particle, and this identifies
the immortal particle as the initial particle. Similarly to the setting
in Section \ref{sec-BRWIP} for IIBRW, we again have a unique infinite line of decent,
which `shakes off' mass at constant rate. This mass which is produced along
the infinite line of decent behaves as a usual super-Brownian motion. This
description can be seen as an alternative construction of \icsbm\!\!.

For \icsbm\!\!,
the immortal particle gives a powerful pictorial description of the moment measures.
Indeed, for $x_1$ in the support of \icsbm at time $s_1$, we can follow back its path
until it hits the path of the immortal particle. Both the path until it hits the
immortal particle's path, as well as the path of the immortal particle before they meet
are Brownian motion paths, and the motion of these two paths will be independent.
According to \refeq{MM1}, the meeting time is uniform
on $[0,s_1]$, which explains the factor $s_1$ in \refeq{MM1}. Thus, the moment measure
$M^{\smallsup{1}}$ can be represented as the union of an infinite path (corresponding to
the immortal particle), together with a path of length $s_1$ which `hooks up' to
the infinite path at a time which is uniform on $[0,s_1]$. We can iterate
this procedure for $x_2$ in the support
of \icsbm at time $s_2$, and this will add a second path from $(x_2, s_2)$
which will `hook up' to one of the two paths present for the moment measure
$M^{\smallsup{1}}$. Iterating this procedure for $M^{\smallsup{r}}$ with $r$
points in the support of \icsbm will create a tree with a single infinite path,
and $r$ paths iteratively connected to the union of the infinite path
and the previously added paths. This picture will
be useful to describe the relation of other models to \icsbm\!\!.

\subsection{Discussion and notes}
The construction of \icsbm presented above is within the
folklore of the super-processes community. For instance, in \cite{EP90},
conditioning on non-extinction was considered in the context of
general measure-valued diffusions starting from a proper
initial measure. However, these results do not immediately
apply to the canonical measure. See \cite{EP90} and the
references therein. The notion of the canonical measure
of super-Brownian motion was first introduced in \cite{EKR91},
and further studied in \cite{LS95}, where its existence and some
of its properties are proved. See also \cite[Theorem 1.4]{DaP99}.

It is well-known that for $d>4$, the range of
super-Brownian motion is four-dimensional. This was first
proved in \cite{DIP89}, see also \cite[Theorem III.3.9]{Perk02}
and the references therein. One would expect that similar estimates
as in \cite[Theorem III.3.9]{Perk02} are also true for \icsbm\!\!.

As far as we know, this is the first time that it was shown that
the moment measures of branching random walk converge to those
of super-Brownian motion. Many of the ideas for this convergence
are taken from \cite{HS02a}, where a similar approach was taken
for oriented percolation above 4 spatial dimensions.
See Section \ref{sec-iicop} below. In \cite{BCHS99}, there is
a related proof that the $r$-point functions of BRW conditioned
to have total mass equal to $n$ converge to the $r$-point functions
of integrated super-Brownian excursion (ISE).

The fact that the total mass at a given time is equal to a
size-biased exponential distribution can also be understood
by \refeq{Ptseq} and the fact that the total mass at a given
time $t$ under the canonical measure conditioned to survive
at time $t$ is an exponential distribution
(see \cite[Theorem II.7.2(iii)]{Perk02}).

\section{The incipient infinite cluster for spread-out
oriented percolation above 4+1 dimensions}
\label{sec-iicop}

Since super-Brownian motion arises as the scaling limit in various critical models,
one can expect that also the \ICSBM arises as a scaling
limit. Here we will give an example where we can prove that \icsbm indeed is the scaling
limit, namely, for the incipient infinite cluster for spread-out oriented percolation
above $4+1$ dimensions. We first discuss the incipient infinite cluster.

In many models, it is known that at the critical value, clusters are finite.
For example, it is believed that the critical percolation probability equals 0.
For oriented percolation on $\Zd \times \Z_+$, it was shown in
\cite{BG90,GH01} that there is no infinite cluster at the critical point.
The notion of the incipient infinite percolation cluster (IIC) is an attempt
to describe the infinite structure that is emerging but not quite present
at the critical point.  Various aspects of the IIC are discussed in \cite{Aize97}.
We start by defining spread-out oriented percolation.

%
%
%

The spread-out oriented percolation models are defined as follows.
Consider the graph with vertices $\Zd \times {\Zbold}_+$ and
directed bonds $((x,n),(y,n+1))$, for $n \geq 0$ and $x,y \in
\Zd$.
Let $D : \Zd \to [0,1]$ be a fixed function.
Let $p \in [0,\|D\|_\infty^{-1}]$, where $\|\cdot\|_\infty$ denotes the
supremum norm, so that $pD(x)\leq
1$ for all $x$. We associate to each directed bond $((x,n),(y,n+1))$ an
independent random variable taking the value $1$ with probability
$pD(y-x)$ and $0$ with probability $1-pD(y-x)$.   We say a
bond is {\em occupied}\/ when the corresponding random variable
is $1$, and {\em vacant}\/ when the random variable
is $0$.
Given a configuration of occupied bonds, we say that $(x,n)$ is {\em connected
to} $(y,m)$, and write $(x,n) \conn (y,m)$, if there is an oriented path
from $(x,n)$ to $(y,m)$ consisting of occupied bonds, or
if $(x,n)=(y,m)$.
The joint probability distribution of the bond variables
will be denoted $\Pbold^{\smallsup{\rm op}}$, with corresponding expectation denoted
$\Ebold^{\smallsup{\rm op}}$.  Note that $p$ is {\em not}\/ a probability, but rather
equals the expected number of occupied bonds per vertex.
We will always work at the critical percolation threshold, i.e.,
at $p=p_c$, and omit subscripts $p_c$ from the notation.

The function
$D$ will always be
assumed to obey the properties of Assumption~D of \cite{HS01a}.
Assumption~D involves a positive parameter $L$, which serves to spread
out the connections, and which we will take to be large.
The parameterisation has been chosen in such a way that $p_c$
will be asymptotically equal to $1$ as $L \to \infty$.
In particular, Assumption~D requires that
$\sum_{x \in \Zd}D(x)=1$, that $D(x) \leq CL^{-d}$ for all $x$, and
that $C_1L \leq \sigma \leq C_2L$ (recall \refeq{sigdef}).

A simple example is
    \eq
    \lbeq{Ddef}
    D(x) = \begin{cases}
    \frac{1}{(2L+1)^{d}-1} & 0< \|x\|_{\infty} \leq L \\
    0 & \mbox{otherwise},
    \end{cases}
    \en
for which bonds are of the form
$((x,n),(y,n+1))$ with $\|x-y\|_{\infty} \leq L$, and
a bond is occupied with probability $p[(2L+1)^{d}-1]^{-1}$.

Let $\Fcal$ denote the $\sigma$-algebra of events.
A {\em cylinder event}\/ is an event that is determined by the occupation
status of a finite set of bonds.  We denote the algebra of cylinder events by
$\Fcal_0$.  Then $\Fcal$ is the $\sigma$-algebra generated by $\Fcal_0$.
For our first definition of the IIC, we begin by defining
${\mathbb P}_n^{\smallsup{\rm op}}$ by
    \eq
    \lbeq{Pndef}
    {\mathbb P}_n^{\smallsup{\rm op}}(E) = \frac{1}{\tau_n^{\smallsup{\rm op}}}
    \sum_{x\in \Z^d} {\mathbb P}^{\smallsup{\rm op}}(E \cap \{ (0,0)\conn (x,n)\})
    \quad (E \in \Fcal_0),
    \en
where $\tau_n^{\smallsup{\rm op}} = \sum_{x \in \Zd}\tau_n^{\smallsup{\rm op}}(x)$ with
$\tau_n^{\smallsup{\rm op}}(x) = \Pbold^{\smallsup{\rm op}}( (0,0) \conn
(x,n))$.
We then define $\Pbold_{\sinfty}^{\smallsup{\rm op}}$ by setting
    \eq\lbeq{IICdef}
    {\mathbb P}_{\sinfty}^{\smallsup{\rm op}}(E)=\lim_{n\rightarrow \infty} {\mathbb P}_n^{\smallsup{\rm op}}(E)
    \quad (E \in \Fcal_0),
    \en
assuming the limit exists. We now turn to the second construction
of the incipient infinite cluster. For this, let
    \eq
    \lbeq{Sndef}
    S_n = \{(0,0)\conn n\}=\{(0,0)\conn (x,n)\;\; \mbox{for some $x\in \Z^d$}\}
    \en
denote the event that the cluster of the origin survives to time $n$.
Define $\Qbold_n$ by
    \eq
    \lbeq{Qndef}
    \Qbold_n^{\smallsup{\rm op}}(E) = \Pbold^{\smallsup{\rm op}} (E | S_n) \quad (E \in \Fcal_0).
    \en
We then define $\Qbold_{\sinfty}^{\smallsup{\rm op}}$ by setting
    \eq
    \lbeq{IICaltdef}
    \Qbold_{\sinfty}^{\smallsup{\rm op}}(E) = \lim_{n \to \infty} \Qbold_n^{\smallsup{\rm op}}(E)
    \quad (E \in \Fcal_0),
    \en
assuming the limit exists.

The following theorem shows that this definition produces a
probability measure on $\Fcal$ under which the origin is almost
surely connected to infinity. The theorem below is the main result
in \cite{HHS02a}. In its statement, we write
    \eq
    \theta_n^\smallsup{{\rm op}}=\Pbold^{\smallsup{\rm op}}(S_n).
    \lbeq{thetandef}
    \en

    \begin{theorem}
    \label{thm-IIC}
    Let $d+1>4+1$ and $p = p_c$. There is an
    $L_0 = L_0(d)$ such that for $L \geq L_0$,
    the limit in \refeq{IICdef} exists
    for every cylinder event $E \in \Fcal_0$.
    Moreover, $\Pbold_\infty^{\smallsup{\rm op}}$ extends to a
    probability measure on the $\sigma$-algebra
    $\Fcal$, and the origin is almost surely connected to
    infinity under ${\mathbb P}_{\infty}^{\smallsup{\rm op}}$.

    If we further assume that there is a finite positive constant
    $B$ such that
        \eq
        \lbeq{OPass}
        \lim_{n\rightarrow \infty} n\theta_n^\smallsup{{\rm op}}=1/B,
        \en
    then also the limit in \refeq{IICaltdef} exists and
    $\Qbold_{\sinfty}^{\smallsup{\rm op}} =
    \Pbold_{\sinfty}^{\smallsup{\rm op}}$.
    \end{theorem}
In \cite{HJ04}, there is
a third construction of the IIC, where we take $p<p_c$, and define
    \eq
    \lbeq{altdefIIC3}
    \mathbb{Q}_{p}^{\smallsup{\rm op}}(E) = \frac{1}{\chi^\smallsup{{\rm op}}(p)} \sum_{(x,n)\in \Z^d\times \Z_+}
    {\mathbb P}^{\smallsup{\rm op}}_p(E \cap \{ (0,0)\conn (x,n)\}),
    \en
where
    \eq
    \chi^\smallsup{{\rm op}}(p) = \sum_{(x,n)\in \Z^d\times \Z_+}\Pbold^{\smallsup{\rm op}}_p
    ((0,0)\conn (x,n))
    \en
denotes the oriented percolation susceptibility, and we now explicitly use the subscript
$p<p_c$ to indicate the percolation parameter.
In \cite{HJ04}, it is proved that when $p\uparrow p_c$, then $\mathbb{Q}_{p}^{\smallsup{\rm op}}$
converges to $\Pbold_{\sinfty}^{\smallsup{\rm op}}$. This definition works both for
oriented and unoriented percolation (see Section \ref{ss:pe}).

\subsection{Convergence of oriented percolation moment measures}
In \cite{HS02a}, it was shown that the finite dimensional distributions of
the rescaled oriented percolation cluster converge to the ones of
super-Brownian motion. We will review this result here. We first define an
analogue random measure valued-process $X_{n,t}$ on
$\R^{dl}$ by placing mass $(A^2Vn)^{-1}$ at each site at times
$\lfloor nt\rfloor$ in $(v\sigma^2n)^{-1/2}C(0,0)$, i.e.,
for any subset $E$ of $\R^d$,
    \eq
    X_{n,t}(E) = \frac{1}{A^2Vn}\sum_{x\in (v\sigma^2n)^{1/2}E} I[(x,\lfloor nt\rfloor)\in C(0,0)].
    \en
Here $v\sigma^2$ serves as the variance of occupied oriented percolation paths.
The main result in \cite{HS02a} is the following:

\begin{theorem}
\label{thm-convSBMOP}
Let $d+1>4+1$. There exist constants $A,V, v$, $\delta\in (0,1)$ and
an $L_0 = L_0(d)$ such that for $L \geq L_0$,
and for all $r \geq 2$ and $\vec{k} \in \R^{d(r-1)}$,
   \eq
   \hat{\tau}_{\vec{n} }^{\smallsup{\rm op}}(\vec{k}/\sqrt{v\sigma^2 n})
    = A (A^2V)^{r-2} n^{r-2}
    \left[
    \hat{M}_{\vec{n}/n}^{\smallsup{r-1}}(\vec{k})
    + {\cal O}((n_{\smallsup{2}}+1)^{-\delta})
        \right]
    \quad
    (r \geq 2)
    \lbeq{aim3ptOP}
    \en
holds uniformly in $n\geq n_{\smallsup{2}}$.
Consequently, $AVn\Pbold(X_{n,\vec{t}} \in \cdot)$ converges to
$\N_0$ in the sense of convergence of finite-dimensional distributions.
\end{theorem}

For a review of the proof of Theorem \ref{thm-convSBMOP}, see \cite[Section 2]{HS01a}.
The approach can be used more generally. For example, the
approach taken in \cite{HSa04c} for the contact process above 4 dimensions,
or the approach taken in \cite{Holm04} for lattice trees, are based upon the
same ideas. We next discuss the convergence of the
oriented percolation IIC $r$-point functions.

The oriented percolation IIC {\em $r$-point functions} are defined,
for $n_i \geq 0$ and $x_i \in \Zd$,  by
    \eq
    \lbeq{mpoint,nobranches}
    \rho^\smallsup{\rm op}_{n_1,\ldots,n_{r-1}}(x_1,\ldots, x_{r-1})
    = {\Pbold}_{\sinfty}^{\smallsup{\rm op}}((0,0) \conn (x_i,n_i)
    \mbox{ for each }
    i=1,\ldots,r-1).
    \en
We now turn to the scaling limit of the oriented percolation IIC. The main result,
which is a direct consequence of results in \cite{HHS02a} and
\cite{HS02a}, is the following theorem:

    \begin{theorem}
    \label{thm-ISBMOP}
    Let $d+1>4+1$. There is an $L_0 = L_0(d)$ such that for $L \geq L_0$,
    and for all $r \geq 2$, $\vec{t} = (t_1,\ldots,t_{r-1}) \in (0,1]^{r-1}$
    and $\vec{k} \in \R^{d(r-1)}$
    \eq
     \lbeq{ISBMrptOP}
        \frac{1}{(A^2V)^{r-1}}
        \hat{\rho}_{m\vec{t}}^\smallsup{{\rm op}} (\vec{k}/\sqrt{v\sigma^2 m})
        =
        m^{r-1} \hat{M}_{\sinfty,\vec{t}}^\smallsup{r}(\vec{k})[1+\Ocal(m^{-\delta})].
    \en
    Consequently, when time is rescaled by $m$ and space by $\sqrt{v\sigma^2 m}$,
    then the finite-dimensional distributions of ${\Pbold}_{\sinfty}^{\smallsup{\rm op}}$
    converge to those of ${\Pbold}_{\sinfty}$.
    \end{theorem}


\proof We will use Proposition \ref{prop-cond}, for which the main assumption
is proved in \cite[(2.52)]{HS02a}. See Theorem \ref{thm-convSBMOP} above.
\qed

In \cite{HHS02a}, there are more properties of the IIC. For instance, a
version of Theorem \ref{thm-sbExp} is proved there, as well as a result
on the four-dimensional nature that we state now. In order to be able to
state the result, we let
    \eq
    C(0,0)=\{(y,m)\in \Z^d\times \Z_+: (0,0) \conn (y,m)\}
    \en
denote the connected cluster of the origin, and let
    \eq
    M^{\smallsup{\rm op}}(R)= \#\{(y,m) \in C(0,0): |y|\leq R\}
    \en
denote the total number of sites in the cluster of the origin that are
at most a distance $R$ away from the origin, under $\Pbold_{\sinfty}$.
We note that $M^{\smallsup{\rm op}}(R)$ is the equivalent of $M(R)$
defined in \refeq{MRdef}. In \cite{HHS02a}, the random variable
$M^{\smallsup{\rm op}}(R)$ has not been studie, but rather its expected value.

    \begin{theorem}
    \label{thm-HD}
    Let $d+1>4+1$ and $p = p_c$.
    There are $L_0 = L_0(d)$ and $C_i=C_i(L,d)>0$ such that
    for $L \geq L_0$,
    \eq
    \lbeq{Drscaling}
    C_1 R^4 \leq \Eb_{\sinfty}^{\smallsup{\rm op}}[M^{\smallsup{\rm op}}(R)] \leq C_2 R^4.
    \en
    \end{theorem}
\noindent
Theorem \ref{thm-HD} is a sign that the IIC is four-dimensional.

We complete this section by showing that essentially there is a unique path
tending to infinity, meaning that any two infinite paths share bonds under
$\Pb_{\sinfty}^{\smallsup{\rm op}}$. This is the equivalent of the immortal particle for
\icsbm and IIBRW. Of course, for oriented percolation, there will be many small doubly
connected parts or sausages along any path to infinity, so that we cannot expect there to be a
{\it unique} infinite path as for IIBRW in Theorem \ref{thm-immparti}.

Before stating the result, we need some definitions. We say that
the events $\{ (y_1,m_1) \conn  (x_1,n_1)\}$ and $\{ (y_2,m_2) \conn (x_2,n_2) \}$
\emph{occur disjointly}, if there exist bond disjoint occupied paths
connecting $(y_1,m_1)$ to $(x_1,n_1)$ and $(y_2,m_2)$ to $(x_2,n_2)$.
We write $\{(y_1,m_1)\conn n\}$ occurs disjointly from
$\{(y_2,m_2)\conn n\}$ for $n\geq m_1\vee m_2$ when there exist $x_1, x_2$ such that
the events $\{ (y_1,m_1) \conn  (x_1,n)\}$ and $\{ (y_2,m_2) \conn (x_2,n)\}$
occur disjointly. We abbreviate this event by $\{(y_1,m_1)\conn n\}\circ \{(y_2,m_2)\conn n\}$
for $n\in \N\cup \{\infty\}$, where $\{(y_1,m_1)\conn \infty\}\circ \{(y_2,m_2)\conn \infty\}$
is the intersection of $\{(y_1,m_1)\conn n\}\circ \{(y_2,m_2)\conn n\}$ for all $n$.

\begin{theorem}
\label{thm-immpartiOP}
Under $\mathbb {P}_{\sinfty}^{\smallsup{\rm op}}$, for every $m$, the probability that
there exist $y_1, y_2\in \Z^d$ such that $(y_1,m)\conn \infty$ occurs disjointly from
$(y_2,m)\conn \infty$ is zero.
\end{theorem}

\proof We bound, using the BK-inequality,
    \eqalign
    &\mathbb {P}_{\sinfty}^{\smallsup{\rm op}}(\{(y_1,m)\conn \infty\}\circ \{(y_2,m)\conn \infty\})\nn\\
    &\qquad= \lim_{k\rightarrow \infty}
    \mathbb {P}_{\sinfty}^{\smallsup{\rm op}}(\{(y_1,m)\conn k\}\circ \{(y_2,m)\conn k\})\nn\\
    &\qquad =\lim_{k\rightarrow \infty}\lim_{n\rightarrow \infty}
    \mathbb {Q}_{n}^{\smallsup{\rm op}}(\{(y_1,m)\conn k\}\circ \{(y_2,m)\conn k\})\nn\\
    &\qquad =\lim_{k\rightarrow \infty}\lim_{n\rightarrow \infty}\frac{1}{\theta_n}\mathbb {P}^{\smallsup{\rm op}}
    \big((\{(y_1,m)\conn k\}\circ \{(y_2,m)\conn k\})\cap \{(0,0)\conn n\})\nn\\
    &\qquad \leq \lim_{k\rightarrow \infty}\lim_{n\rightarrow \infty} \frac{2\theta_{n}\theta_{k-m}}{\theta_n}
   =0,
    \enalign
where we use the fact that
    \eqalign
    &(\{(y_1,m)\conn k\}\circ \{(y_2,m)\conn k\})\cap \{(0,0)\conn n\}\nn\\
    &\qquad \subseteq
    (\{(y_1,m)\conn k\}\circ \{(0,0)\conn n\})\cup (\{(y_2,m)\conn k\}\circ\{(0,0)\conn n\}).
    \enalign
Since the event that there exist $y_1, y_2$ such that $(y_1,m)\conn \infty$
occurs disjointly from $(y_2,m)\conn \infty$ is a countable union of
events with probability 0, the claim follows.\qed


\section{Conjectured scaling to \icsbm\!\!: Incipient structures}
\label{sec-conjII}
In this section, we describe several models of incipient infinite
structures where one can expect convergence to \icsbm to hold.

\subsection{The incipient infinite cluster for
percolation above 6 dimensions}
\label{ss:pe}
For general background on percolation, see~\cite{Grim99}.
Our models are defined in terms of a function
$D: \Zd \to [0,1]$. Let $p \in [0, \| D \|_\infty^{-1}]$ be a parameter, so that
again $p D(x) \le 1$ for all $x$. We declare a bond $\{u,v\}$ to be {\em occupied}
with probability $p D(v - u)$ and {\em vacant} with probability
$1 - p D(v - u)$. The occupation status of all bonds are
independent random variables. For the nearest-neighbor model,
we take $D(x)=1/(2d)$ for
all $x$ with $|x|=1$, so that each bond is occupied with probability
$p/(2d)$. For the spread-out model, we assume that the conditions in
\cite[Definition 1.1]{HHS01a} are satisfied.
The function in~\refeq{Ddef} does obey the assumptions.

The law of the configuration of occupied bonds (at the critical
percolation threshold)
is denoted by $\Pb^{\smallsup{\rm pe}}$ with corresponding expectation denoted by
$\Eb^{\smallsup{\rm pe}}$. Given a configuration we say that $x$ is connected to $y$,
and write $x \cnctd y$, if there is a path of occupied bonds from $x$ to $y$
(or if $x = y$).

Let $\Fcal$ denote the $\sigma$-algebra of events. A {\em cylinder event}
is an event given by conditions on the states of finitely many bonds only.
We denote the algebra of cylinder events by $\Fcal_0$. We define
\eq
\lbeq{def-P_x}
  \Pb_x^{\smallsup{\rm pe}}(F) = \Pb^{\smallsup{\rm pe}}(F|0 \cnctd x)
           = \frac{1}{\tau^{\smallsup{\rm pe}}(x)} \Pb^{\smallsup{\rm pe}}(F,\, 0 \cnctd x),
         \quad F \in \Fcal,
\en
where $\tau^{\smallsup{\rm pe}}(x) = \Pb^{\smallsup{\rm pe}}(0 \cnctd x)$. The main result in
\cite{HJ04} is the following theorem:

\begin{theorem}
\label{thm:IIC-lim}
Let $d > 6$ and $p = p_c$. There is an $L_0 = L_0(d)$ such that for
$L \ge L_0$ in the spread-out model, the limit
\eq
\lbeq{IIC-lim}
  \Pb^{\smallsup{\rm pe}}_{\sinfty}(F) = \lim_{|x| \to \infty}\Pb_x^{\smallsup{\rm pe}}(F)
\en
exists for any cylinder event $F$. Also, $\Pb^{\smallsup{\rm pe}}_{\sinfty}$ extends
uniquely from $\Fcal_0$ to a probability measure on $\Fcal$.
\end{theorem}


Theorem \ref{thm:IIC-lim} is similar to the existence statement of the
IIC for spread-out oriented percolation above $4+1$ dimensions in
Theorem \ref{thm-IIC}. Moreover, the definition in \refeq{altdefIIC3} is
also proved to exist in \cite{HJ04}, and to give the same result.
In \cite{HJ04} some properties of $\Pb_{\sss \infty}^{\smallsup{\rm pe}}$
were proved, which are the natural equivalents of Theorems
\ref{thm-HD}--\ref{thm-immpartiOP}.

We next turn to the conjecture linking the unoriented percolation IIC to
\icsbm\!\!. Of course, there is no explicit
time variable in unoriented percolation, so will
will introduce a natural candidate for a time variable. Define
$\SP(x,y)$ to be the shortest path along occupied bonds between
$x$ and $y$, and let $|\SP(x,y)|$ be the number of bonds in this
shortest path. When $x$ and $y$ are not connected, then we set
$|\SP(x,y)|=\infty$. Let $\SP(x)=\SP(0,x)$. We then think of
$|\SP(x)|$ as being a time variable analogous
to the time variable $n$ in oriented percolation. Define
    \eq
    \tau^{\smallsup{\rm pe}}_{\vec{n}}(\vec{x})
    = {\Pbold}^{\smallsup{\rm pe}}(|\SP(x_j)|=n_j
    \mbox{ for each }
    j=1,\ldots,r-1)
    \en
and
    \eq
    \rho^{\smallsup{\rm pe}}_{\vec{n}}(\vec{x})
    = {\Pbold}_{\sinfty}^{\smallsup{\rm pe}}(|\SP(x_j)|=n_j
    \mbox{ for each }
    j=1,\ldots,r-1).
    \en
Then $\tau^{\smallsup{\rm pe}}_{\vec{n}}(\vec{x})$ is analogous to
the oriented percolation probability
$\tau^{\smallsup{\rm op}}_{\vec{n}}(\vec{x})$ of \refeq{mpoint,nobranches},
while $\rho^{\smallsup{\rm pe}}_{\vec{n}}(\vec{x})$ is analogous to
$\rho^{\smallsup{\rm op}}_{\vec{n}}(\vec{x})$. Then we conjecture that
$\rho^{\smallsup{\rm pe}}_{\vec{n}}$ converges to the moment measures of
\icsbm\!\!:

    \begin{conj}
    \label{conj-ISBMPE}
    Let $d>6$. For all $r \geq 2$, $\vec{t} = (t_1,\ldots,t_{r-1})
    \in \R^{r-1}$ and $\vec{k} \in \R^{d(r-1)}$, there exist constants $A,V,v$
    and $\delta\in (0,1)$ such that \refeq{rhoscal} holds for $\rho^{\smallsup{\rm pe}}$.
    \end{conj}
In order to prove Conjecture \ref{conj-ISBMPE}, the key step is to prove a version of
\refeq{aim3pt} for unoriented percolation, as conjectured in \cite[Section 1.3.3]{HS02a}.

\subsection{The contact process incipient infinite cluster
above 4 dimensions}
\label{sec-iiccp}
For a general introduction to the contact process, see \cite{Ligg99}.
We define the spread-out contact process as follows.
Let $\bC_t\subset\Z^d$ be the set of infected individuals
at time $t\in\R_+$, and let $\bC_0=\{0\}.$ An infected
site $x$ recovers in a small time interval $[t,t+\vep]$
with probability $\vep+o(\vep)$ independently of $t$,
where $o(\vep)$ is a function that satisfies
$\lim_{\vep\downarrow 0}o(\vep)/\vep=0$. In other words, $x\in\bC_t$
recovers with rate 1.  A healthy site $x$ gets infected, depending
on the status of its neighboring sites, with rate $\lamb\sum_{y\in\bC_t}D(x-y)$,
where $\lamb\geq0$ is the infection rate.  We denote
the associated probability measure by $\mP^\lamb$.  We will assume that
the function $D:\Zd\mapsto[0,1]$ is a probability distribution which
satisfies the assumptions in Section \ref{sec-iicop}.

We will always investigate the contact process at the critical value
$\lambc$ for the sufficiently spread-out contact process above 4 dimensions.
In \cite{HSa04a, HSa04c}, the contact process is investigated and
the goal is to prove a version of \refeq{aim3pt} for the contact process.
For this, its close analogy to oriented percolation is essentially used.
We now explain this connection.

The contact process can be constructed using a graphical
representation.  We consider $\Zd\times\R_+$ as
space-time.  Along each time line $\{x\}\times\R_+$, we place
points according to a Poisson process with intensity 1,
independently of the other time lines.  For each ordered pair of
distinct time lines from $\{x\}\times\R_+$ to $\{y\}\times\R_+$,
we place directed bonds $((x,\,t),(y,\,t))$, $t\geq0$, according to
a Poisson process with intensity $\lamb\,D(y-x)$,
independently of the other Poisson processes.  A site $(x,s)$ is
said to be {\it connected to} $(y,t)$ if either $(x,s)= (y,t)$ or
there is a non-zero path in $\Zd\times\R_+$ from $(x,s)$ to
$(y,t)$ using the Poisson bonds and time line segments traversed
in the increasing time direction without traversing the Poisson
points.  The law of $\bC_t$ defined above is
equal to that of $\{x\in \Zd:(0,0)$ is connected to
$(x,\,t)\}$.

Inspired by this percolation structure in space-time
and following \cite{s01}, we consider
the following oriented percolation process in $\Zd\times\vep\Z_+$ with
$\vep \in(0,1]$ being a discretization parameter.  A directed pair
$b=((x,t),(y,t+\vep))$ of sites in $\Zd\times\vep\Z_+$ is called a
{\it bond}.  Each bond is either {\it
occupied} or {\it vacant} independently of the other bonds, and a
bond $b=((x,t),(y,t+\vep))$ is occupied with probability
    \eq
    \lbeq{bprob} p_\vep(y-x)=\begin{cases}
    1-\vep,&\mbox{if }x=y,\\
    \lamb\vep\,D(y-x),&\mbox{otherwise},
    \end{cases}
    \en
provided that $\sup_x p_{\vep}(x)\leq1$.  We denote the
associated probability measure by $\mP_\vep^\lamb$.  It is proved in
\cite{BG91} that $\mP_\vep^\lamb$ weakly converges to $\mP^\lamb$ as
$\vep\downarrow0$.

Existence of the contact process IIC has not yet been established.
The proof in \cite{HHS02a} applies to the discretized contact
process, and therefore, the only thing left to do is to take the limit
$\vep\downarrow 0$. The continuum limit results in \cite{HSa04a}
can hopefully show that this continuum limit exists. Once the existence
of the contact process IIC has been established, convergence of the
finite-dimensional distributions towards the moments measures of
\icsbm will follow from the results in \cite{HSa04a}.

\subsection{Incipient infinite lattice trees above 8 dimensions}
There are many connections between lattice trees and super-Brownian motion.
In \cite{DS97, DS98}, it was shown that the $r$-point functions of lattice trees
of fixed size, converge to those of ISE. The statements are complete when dealing
with the $r$-point functions where the number of steps between $0$ and $x$ along the
tree is not fixed, and there are partial results when this number is fixed and
scales with the size of the lattice tree.

There is current progress in
understanding the connection to SBM \cite{Holm04}, when the set-up is somewhat
different. Let us introduce some notation. A {\it lattice tree} is a tree embedded
in $\Z^d$ containing no cycles. We give uniform weight to lattice trees with a
fixed number of bonds, and assume that the bonds are either nearest-neighbour,
or spread-out (as in \refeq{Ddef}). In general, the number of lattice trees
of fixed size grows exponentially with the size. Denote by $\tau^{\smallsup{\rm lt}}(N)$
the total number of lattice trees of size $N$ containing $0$.
Then, we know that
    \eq
    \lim_{N\rightarrow \infty} \tau^{\smallsup{\rm lt}}(N)^{1/N} =\lambda\in (0,\infty).
    \en
We define
    \eq
    \lbeq{rptLT}
    \tau^{\smallsup{\rm lt}}_{\vec{n}}(\vec{x}) =
    \sum_{N=1}^{\infty} \tau^{\smallsup{\rm lt}}_{\vec{n}}(\vec{x};N)\lambda^{-N}
    \en
to be the $r$-point function for lattice trees, where $\tau^{\smallsup{\rm lt}}_{\vec{n}}(\vec{x};N)$
is the number of lattice trees of size $N$ such that $|\SP(x_j)|=n_j$ for all
$j=1, \ldots, r-1$. Of course, existence of the
sum in \refeq{rptLT} is a non-trivial result, and follows from \cite{HHS01a} for the spread-out model.
The main work in \cite{Holm04} involves the proof that $\tau^{\smallsup{\rm lt}}_{\vec{n}}(\vec{x})$
scales to the $r$-point functions of SBM, and therefore satisfies the main assumption
in \refeq{aim3pt}. If this is completed, it is natural to conjecture that there exists
an infinite tree measure, and that this infinite tree measure has \icsbm as a scaling limit.
The simplest way to obtain the infinite tree measure is to take the limit
    \eq
    \lbeq{IT}
    \Pbold^{\smallsup{\rm lt}}_{\sinfty}(E)=\lim_{n\rightarrow \infty}
    \frac{\sum_x\tau^{\smallsup{\rm lt}}_{n}(x,E)}{\sum_x\tau^{\smallsup{\rm lt}}_{n}(x)},
    \en
where, for an event $E$,
    \eq
    \tau^{\smallsup{\rm lt}}_{n}(x,E)=\sum_{N=1}^{\infty}\tau^{\smallsup{\rm lt}}_{n}(x,E;N)\lambda^{-N},
    \en
and $\tau^{\smallsup{\rm lt}}_{n}(x,E;N)$ is the number of lattice trees of size $N$
such that $|\SP(x)|=n$ and such that $E$ is satisfied. $\Pbold^{\smallsup{\rm lt}}_{\sinfty}$ is
a version of $\Pbold_{t}$ in \refeq{Ptdef}. We next state a version of $\Qbold_{t}$ in \refeq{Qtdef}.
For this, we let $\tau^{\smallsup{\rm lt}}_{n}(E;N)$ denote the number of lattice trees such that
there exists an $x$ with $|\SP(x)|\geq n$ satisfying $E$, and
    \eq
    \tau^{\smallsup{\rm lt}}_{n}(E) = \sum_{N=1}^{\infty}
    \tau^{\smallsup{\rm lt}}_{n}(E;N) \lambda^{-N}.
    \en
Finally, let
    \eq
    \theta^{\smallsup{\rm lt}}_{n}=\tau^{\smallsup{\rm lt}}_{n}(\Omega),
    \en
where $\Omega$ is the whole probability space. Then we define
    \eq
    \lbeq{IT2}
    \Qbold^{\smallsup{\rm lt}}_{\sinfty}(E)=\lim_{n\rightarrow \infty}
    \frac{\tau^{\smallsup{\rm lt}}_{n}(E)}{\theta^{\smallsup{\rm lt}}_{n}},
    \en
assuming the limit exists.

It should be possible to use
the lace expansion to prove that the limit in \refeq{IT} is well-defined, but this has not
yet been done. The limit in \refeq{IT2} will be much more involved,
since for this, one needs to understand the lattice tree survival probability
$\theta^{\smallsup{\rm lt}}_{n}$. It is
natural to conjecture that
    \eq
    \rho_{\vec m}^\smallsup{{\rm lt}} (\vec{x})=\Pbold^{\smallsup{\rm lt}}_{\sinfty}
    (|\SP(x_j)|=m_j\forall j=1,\ldots, r-1)
    \en
scales to the $r$-point function of \icsbm for $d>8$. That is the content of the next conjecture:

    \begin{conj}
    \label{conj-ISBMLT}
    Let $d>8$. For all $r \geq 2$, $\vec{t} = (t_1,\ldots,t_{r-1})
    \in \R^{r-1}$ and $\vec{k} \in \R^{d(r-1)}$, there exist constants $A,V,v$
    and $\delta\in (0,1)$ such that \refeq{rhoscal} holds for $\rho^{\smallsup{\rm lt}}$.
    \end{conj}

\section{Conjectured scaling to \icsbm\!\!: Infinite structures}
\label{sec-conjI}
So far, we have given a number of conjectures linking \icsbm to
incipient infinite structures. We end this paper with two examples where the
structures are infinite.

\subsection{Invasion percolation above 6 dimensions}

We introduce the model for invasion percolation. For simplicity,
we only define the model for a uniform step distribution $D$,
such as the nearest-neighbour case or the case in \refeq{Ddef}.
The bonds in these models are ${\mathbb B}=\{b=(u,v): D(u-v)>0\}$. We let
$\{\omega(b)\}_{b\in \Bbold}$ be a collection of i.i.d.\ uniform random
variables. Given a random configuration $\omega$, we define a random
increasing sequence of subgraphs $G_0, G_1, \ldots$ as follows.
We let $G_0$ be the graph with no edges, and the single vertex $0$. We
let $G_{i+1}=G_i \cup\{b_{i+1}\}$, where the edge $b_{i+1}$ is obtained
by taking the $b \notin G_i$ with minimal $\omega(b)$ and such that
$b$ has an end vertex in $G_i$. The invaded region is
${\cal S}=\cup_{i=0}^{\infty} G_i$.
The law of the configuration of bonds in the invaded region
is denoted by $\Pb^{\smallsup{\rm ip}}$ with corresponding expectation denoted by
$\Eb^{\smallsup{\rm ip}}$.

It is well-known that the asymptotic behaviour of invasion
percolation is closely related to the incipient cluster.
The heuristic behind this is that
$\limsup_{i \to \infty} \omega(b_{i}) = p_c$,
which is the critical percolation threshold in the model \cite{CCN85}.
In other words, asymptotically the invasion process only accepts
values from critical clusters. As mentioned earlier, critical
clusters in $d > 6$ are four-dimensional, which leads to the
well-known conjecture \cite{NS96} that
$\mathbb{P}( \text{$y$ is invaded} ) \asymp |y|^{-(d-4)}$ when $d > 6$.
This conjecture is supported by results in \cite{HJ04}, where the
conjecture that
    \eq
    \Pbold (\text{$y$ is invaded})
    \asymp \Pbold_{\sss \infty} (0 \conn y)
    \asymp |y|^{-(d -4)}
    \lbeq{con2}
    \en
is explained in some detail. We can stretch this conjecture much further, and
conjecture that the scaling limit of invasion percolation above 6 dimensions
is \icsbm\!\!. For this, we define
    \eq
    \rho^{\smallsup{\rm ip}}_{\vec{n}}(\vec{x})
    = {\Pbold}^{\smallsup{\rm ip}}(|\SP(x_j)|=n_j
    \mbox{ for each }
    j=1,\ldots,r-1),
    \en
where now $|\SP(x)|$ is the minimal number of bonds in the invaded region
along paths from $0$ to $x$. Then we conjecture that
$\rho^{\smallsup{\rm ip}}_{\vec{n}}$ converges to the moment measures of
\icsbm\!\!:

    \begin{conj}
    \label{conj-ISBMIP}
    Let $d>6$. For all $r \geq 2$, $\vec{t} = (t_1,\ldots,t_{r-1})
    \in \R^{r-1}$ and $\vec{k} \in \R^{d(r-1)}$, there exist constants $A,V,v$
    and $\delta\in (0,1)$ such that \refeq{rhoscal} holds for $\rho^{\smallsup{\rm ip}}$.
    \end{conj}

   Conjecture \ref{conj-ISBMIP} is quite hard to prove, as the relation between
   invasion percolation and unoriented percolation is not very direct.
   Maybe it would be simpler to investigate the following variant of invasion
   percolation, where we let $\{\omega(b)\}_{b\in \Bbold}$ be a collection of
   i.i.d.\ uniform random variables on $[0,p_c]$ with probability $p_c$ and are equal
   to $\infty$ with probability $1-p_c$. In this case, instead of picking
   the smallest weight that is larger than $p_c$ when none below $p_c$ is available,
   we simply pick each of the boundary bonds of $G_i$ with equal probability.

\subsection{Uniform spanning forest above 4 dimensions}
The uniform spanning forest (USF) can be obtained as the weak limit of ordinary wired
spanning trees on a large cube when the size of the cube tends to infinity. As
it turns out, for $d\leq 4$, the USF consists of a single tree, while for $d>4$
it consists of multiple trees. See \cite{BKPS02} and the references therein.
In \cite{BKPS02}, there is a wealth of properties of USF's. For instance, a.s., the
maximum over $x$ and $y$ of the number of edges outside the USF in a path from $x$
to $y$ equals $\lfloor \frac{d-1}{4}\rfloor$. Also, in \cite{BKPS02}, it is shown that
the USF has stochastic dimension 4, which is a version of the statement that the trees
that the USF consists of are four-dimensional.

Single uniform spanning trees, such as the spanning trees containing the origin,
are natural candidates for convergence to \icsbm\!\!. To explain this in more detail,
we need some notation. For $x\in \Z^d$, we let $T(x)$ be the (infinite) tree that contains $x$.
Also, for $x\in T(y)$, we let ${\rm SP}(y,x)$ denote the path in the tree $T(x)=T(y)$ that goes
from $x$ to $y$, and we let ${\rm SP}(y,x)|$ denote the number of bonds in this path
and write ${\rm SP}(x)={\rm SP}(0,x)$. Then, we define the $r$-point functions to be
    \eq
    \rho_{\vec n}^{\smallsup{\rm st}}(\vec x)
    =\mathbb{P}(x_i\in T(0), |{\rm SP}(x_i)|=n_i \forall i=1, \ldots, r-1).
    \en

``Wilson's method rooted at infinity'' \cite{Wils96} can be used to generate the
shortest path tree between any number of points using loop-erased random walk,
and thus, allows us to give a probabilistic representation for the event that
$x_i\in T(0)$ and $|{\rm SP}(x_i)|=n_i$ for all $i=1, \ldots, r-1$.
Wilson's method rooted at infinity works on any transient graph. We start by generating
an infinite simple random walk from the origin, and loop-erase it. Call the result
$F_1$. Then, we start an simple random walk from $x_1$, and stop it when it hits
$F_1$ (it is possible that it does not hit $F_1$ at all). After this, we loop-erase it,
and call the union of the two loop-erased paths $F_2$. We can iterate this procedure.
Denote by $F_{k-1}$ the union of the loop-eared paths from $0$, $x_1, \ldots, x_{k-1}$.
Then we start a simple random walk from $x_k$ until it hits $F_{k-1}$, and subsequently
loop-erase it, giving a self-avoiding path $\gamma_k$. Denote $F_k=F_{k-1}\cup \gamma_k,$
and repeat the above procedure. To obtain the USF, we will have to go through all
points of the graph $\Z^d$. However, the order in which the points are chosen is irrelevant,
and thus, for the $r$-point function, it is convenient to start with $0, x_1, \ldots, x_{r-1}.$
Clearly, when $x_i\in T(0)$ for all $i$, the result is the shortest path tree contained
in $T(0)$ containing $x_1, \ldots, x_{r-1}$. Thus, $\rho_{\vec n}^{\smallsup{\rm st}}(\vec x)$
equals the probability that $x_i\in T(0)$, and moreover, that the distance in $T(0)$ between
$0$ and $x_i$ equals $n_i$.

Wilson's construction shows that the behaviour of the USF is intimately related to
loop-erased random walks (LERW). A lot is known about LERW, especially in dimensions
$d\geq 4$. For $d>4$, LERW behaves diffusively, and the rescaled path converges to Brownian
motion. See \cite[Chapter 7]{Lawl91} and the references therein. Therefore, we can think
of the tree containing the origin and the points $x_1, \ldots, x_r$ as built up from an
infinite path which scales to a Brownian motion (the LERW starting at the origin), and
$r$-paths which also scale to Brownian motion and that are iteratively added to
the infinite path. This picture of the $r$-point functions for the USF agrees with the
picture of the $r$-point functions of \icsbm using the immortal particle.
This leads us to the following conjecture:

    \begin{conj}
    \label{conj-ISBMST}
    Let $d>4$. For all $r \geq 2$, $\vec{t} = (t_1,\ldots,t_{r-1})
    \in \R^{r-1}$ and $\vec{k} \in \R^{d(r-1)}$, there exist constants $A,V,v$
    and $\delta\in (0,1)$ such that \refeq{rhoscal} holds for $\rho^{\smallsup{\rm st}}$.
    \end{conj}


\subsection{Discussion and notes for Sections \ref{sec-iicop}--\ref{sec-conjI}}
We have presented a number of examples where the scaling limit is proven or
conjectured to be \icsbm\!\!. Therefore, \icsbm is a natural and
robust object that arises as a universal limit in a variety of
models. We have focussed on convergence of the $r$-point functions
to the moment measure of \icsbm\!\!. It would be of interest to
prove stronger versions of convergence and to prove tightness.
Tightness has proved to be difficult in all the models we have
described in this paper.

The first results showing that spread-out oriented percolation
above 4 spatial dimensions is Gaussian can be found in \cite{NY93, NY95},
where the triangle condition is verified, and the two-point function
is studied. The relation between oriented percolation and SBM in
\cite{HS02a} was the greatest source of inspiration for this paper.
Many more properties for oriented percolation, and its close brother,
the contact process (see Section \ref{sec-iiccp} below) are known.
For example, in \cite{s02}, hyperscaling inequalities are derived
for these two models. From these hyperscaling inequalities, it follows
that the mean-field critical exponents are restricted to $d\geq 4$.
This identifies 4 as the upper critical dimension for oriented percolation.
In \cite{HHS02a}, many related properties are proved for critical
spread-out oriented percolation. For example, it is shown that
the mass at time $m$, properly rescaled as in Theorem \ref{thm-NmBRW},
converges to a size-biased exponential random variable, and that,
conditioned to be alive at time $m$, $N_m/m$ weakly converges to an
exponential random variable. Theorem \ref{thm-immpartiOP} is not
proved in \cite{HHS02a}, even though all the tools were available
at that point.

The assumption in \refeq{OPass} is under investigation in \cite{HHS04a, HHS04b},
where we investigate the critical survival probability, using the lace expansion.


For incipient infinite structures, we have a skeleton of a proof
for the convergence towards \icsbm\!\!, by using convergence
to SBM and Theorem \ref{prop-cond}. However, for infinite structures, this
approach cannot be followed. Therefore, one would have to work with
the $r$-point functions of the infinite structures directly, and
prove convergence by investigating their scaling. It would be interesting,
but probably quite difficult, to derive this scaling for one of the
two examples.

For loop-erased random walk, it is known that also in the upper
critical dimension, convergence towards Brownian motion holds,
with logarithmic corrections. It would be of interest, but
probably quite difficult, to extend this result to convergence
towards \icsbm\!\!.

In \cite{CDP00}, it is shown that the rescaled finite-range voter model
converges to SBM for $d\geq 2$. Also, several related results, where
local mean-field limits are taken, or for critical (continuous time)
branching random walks, are considered. The proofs of these results are stronger
than the ones for oriented percolation, since also tightness is proved.
The reason that this problem is simpler is the fact that the dual process
is coalescing random walks, and this is a simpler process. Thus, martingale
methods can be used to prove the convergence to SBM. It would be of interest to
investigate the link with the canonical measure further. Also, it can be expected
that the voter model conditioned on non-extiction, and where we start with a single
person having different opinion from all others, converges to \icsbm\!\!.
Similar methods as in \cite{CDP00} are used
in \cite{DP99}, where it was shown that the contact process converges to
SBM when the range of the process grows with time. This mean-field limit with
growing ranges makes the problem simpler than in the case where the range is large,
but {\it fixed}, and this is the reason that SBM already appears as the limit
when $d\geq 2$, rather than for $d>4$ as in Section \ref{sec-iiccp}.

We close this discussion with an example of a tree in $\Z^d$ that
is conjectured {\it not} to scale to \icsbm\!\!\!. For this,
give each bond in $\Z^d$ a uniform weight in $[0,1]$.
The {\it minimal spanning forest} (MSF) is the subgraph
on $\Z^d$ where from each cycle we remove the edge with
maximal weight. Then, it is conjectured that the stochastic
dimension of the resulting structure equals 8 above 8 dimensions
(see \cite[Conjecture 6.7]{BKPS02}. Since \icsbm is 4-dimensional,
it cannot be expected that \icsbm can arise as the scaling limit
of MSF. It would be of interest to (even heuristically) identify
the scaling limit of MSF.

\section*{Acknowledgements}
This work was supported in part by Netherlands Organisation for
Scientific Research (NWO). I thank Gordon Slade for introducing me to the
wonderful world of scaling limits and super-processes, and for being
able to explore this world together with him over the past years.
Also, I thank Gordon for several other useful suggestions to
improve the presentation in this paper.
I also thank Takashi Hara, Frank den Hollander, Antal J\'arai and Akira Sakai, who
have been my collaborators on research projects that have led to this
article, and Frank Redig for stimulating discussions concerning the USF.
I thank Mark Holmes and Akira Sakai for useful remarks on early
versions of this paper. Finally, I thank Achim Klenke,
Jean-Fran\c{c}ois Le Gall and Ed Perkins for clarifying
some of the subtleties of super-processes.

\end{document}